\pdfoutput=1
\RequirePackage{ifpdf}
\ifpdf 
\documentclass[pdftex]{sigma}
\else
\documentclass{sigma}
\fi

\usepackage{bm}
\usepackage{booktabs}
\usepackage{multirow}
\usepackage{stmaryrd}
\usepackage{tensor}
\usepackage[vcentermath]{youngtab}

\numberwithin{equation}{section}

\newtheorem{Theorem}{Theorem}[section]
\newtheorem{Corollary}[Theorem]{Corollary}
\newtheorem{Lemma}[Theorem]{Lemma}
\newtheorem{Proposition}[Theorem]{Proposition}
\newtheorem{Problem}[Theorem]{Problem}
\newtheorem{Question}[Theorem]{Question}
\newtheorem{Conjecture}[Theorem]{Conjecture}
 { \theoremstyle{definition}
\newtheorem{Definition}[Theorem]{Definition}
\newtheorem{Remark}[Theorem]{Remark} }

\renewcommand{\ge}{\geqslant}
\renewcommand{\le}{\leqslant}
\newcommand{\coloneq}{:=}
\newcommand{\eqcolon}{=:}
\newcommand{\isom}{\cong}
\newcommand{\transpose}{t}
\newcommand{\Hodge}{\ast}
\newcommand{\adjoint}{\star}
\newcommand{\Id}{I}

\DeclareMathOperator{\tr}{tr}
\DeclareMathOperator{\Tr}{Tr}
\DeclareMathOperator{\Sym}{Sym}
\DeclareMathOperator{\End}{End}
\DeclareMathOperator{\GLG}{GL}
\DeclareMathOperator{\OG}{O}
\DeclareMathOperator{\SOG}{SO}
\DeclareMathOperator{\SpinG}{Spin}
\DeclareMathOperator{\Adj}{Adj}
\DeclareMathOperator{\AV}{\mathcal V}
\DeclareMathOperator{\PV}{\mathbb P\AV}

\newcommand{\aone}{\scalebox{\scaleboxfactor}{\ensuremath{a_1}}}
\newcommand{\atwo}{\scalebox{\scaleboxfactor}{\ensuremath{a_2}}}
\newcommand{\bone}{\scalebox{\scaleboxfactor}{\ensuremath{b_1}}}
\newcommand{\btwo}{\scalebox{\scaleboxfactor}{\ensuremath{b_2}}}
\newcommand{\cone}{\scalebox{\scaleboxfactor}{\ensuremath{c_1}}}
\newcommand{\ctwo}{\scalebox{\scaleboxfactor}{\ensuremath{c_2}}}
\newcommand{\done}{\scalebox{\scaleboxfactor}{\ensuremath{d_1}}}
\newcommand{\dtwo}{\scalebox{\scaleboxfactor}{\ensuremath{d_2}}}

\newcommand{\scaleboxfactor}{0.70710678118654746}

\begin{document}

\newcommand{\arXivNumber}{1205.6227}


\renewcommand{\PaperNumber}{080}

\FirstPageHeading

\ShortArticleName{The Variety of Integrable Killing Tensors on the 3-Sphere}

\ArticleName{The Variety of Integrable Killing Tensors\\ on the 3-Sphere}

\Author{Konrad SCH\"OBEL}

\AuthorNameForHeading{K.~Sch\"obel}

\Address{Institut f\"ur Mathematik,
	Fakult\"at f\"ur Mathematik und Informatik,\\
	Friedrich-Schiller-Universit\"at Jena,
	07737 Jena,
	Germany}

\Email{\href{mailto:konrad.schoebel@uni-jena.de}{konrad.schoebel@uni-jena.de}}

\ArticleDates{Received November 14, 2013, in f\/inal form July 15, 2014; Published online July 29, 2014}

\Abstract{Integrable Killing tensors are used to classify orthogonal coordinates in which the classical Hamilton--Jacobi equation can
	be solved by a separation of variables.  We completely solve the Nijenhuis integrability conditions for Killing tensors on
	the sphere $S^3$ and give a set of isometry invariants for the integrability of a Killing tensor.  We describe explicitly
	the space of solutions as well as its quotient under isometries as projective varieties and interpret their
	algebro-geometric properties in terms of Killing tensors.  Furthermore, we identify all St\"ackel systems in these varieties.
	This allows us to recover the known list of separation coordinates on $S^3$ in a simple and purely algebraic way.  In
	particular, we prove that their moduli space is homeomorphic to the associahedron $K_4$.}


\Classification{53A60; 14H10; 14M12}

{\small \tableofcontents}

\section{Introduction}

\subsection{The problem}

\looseness=-1
Given a partial dif\/ferential equation, it is a natural problem to seek for coordinate systems in which this equation can be
solved by a separation of variables and to classify all such coordinate systems.  It is not surprising, that for classical
equations this problem has a very long history and is marked by the contributions of a number of prominent mathematicians,
beginning with Gabriel Lam\'e in the f\/irst half of of the 19th century \cite{Lame}.  For example, consider the Hamilton--Jacobi
equation
\begin{gather}
	\label{eq:Hamilton-Jacobi}
	\frac12g^{ij}\frac{\partial S}{\partial x^i}\frac{\partial S}{\partial x^j}+V(x)=E
\end{gather}
on a Riemannian manifold.  We say that it separates (additively) in a system of coordinates $x_1,\ldots,x_n$ if it admits a
solution of the form
\[
	S(x_1,\ldots,x_n;\underline c)=\sum_{i=1}^nS_i(x_i;\underline c),
	\qquad
	\det\left(\frac{\partial^2S}{\partial x^i\partial c_j}\right)\not=0,
\]
depending on $n$ constants $\underline c=(c_1,\ldots,c_n)$.  If this is the case, the initial partial dif\/ferential equation
decouples into $n$ ordinary dif\/ferential equations.  We call the corresponding coordinate system \emph{separation coordinates}.

It turns out that the Hamilton--Jacobi equation in the form \eqref{eq:Hamilton-Jacobi} separates in a given ortho\-go\-nal coordinate
system if and only if it also separates for $V\equiv0$ and if the potential $V$ satisf\/ies a certain compatibility condition.  In
fact, the Hamilton--Jacobi equation with $V\equiv0$ plays a~key role for separation of variables in general.  The reason is that
for orthogonal coordinates the (additive) separation of this equation is a necessary condition for the (multiplicative)
separation of other classical equations such as the Laplace equation, the Helmholtz equation or the Schr\"odinger equation.
Suf\/f\/icient conditions can be given in the form of compatibility conditions for the curvature tensor and for the potential.

The classif\/ication of orthogonal separation coordinates for the Hamilton--Jacobi equation thus leads to a variety of ordinary
dif\/ferential equations.  These are well-known equations that def\/ine special functions appearing all over in mathematics and
physics~-- such as Bessel functions, Le\-gendre polynomials, spherical harmonics or Mathieu functions, to name just a few.  In
particular, these functions serve as bases for series expansions in the explicit solution of boundary value problems on domains
bounded by coordinate hypersurfaces.

Separation coordinates are intimately related to the existence of rank two Killing tensors that satisfy a certain integrability
condition.
\begin{Definition}
	A symmetric tensor $K_{\alpha\beta}$ on a Riemannian manifold is a {\it Killing tensor} if and only if
	\begin{gather}
		\label{eq:Killing}
		\nabla_{(\gamma}K_{\alpha\beta)}=0		 ,
	\end{gather}
	where $\nabla$ is the Levi-Civita connection of the metric $g$ and the round parenthesis denote complete symmetrisation of
	the enclosed indices.  Depending on the context, we will consider a Killing tensor either as a symmetric bilinear form
	$K_{\alpha\beta}$ or alternatively as a symmetric endomorphism $K\indices{^\alpha_\beta}=g^{\alpha\gamma}K_{\gamma\beta}$.
\end{Definition}

\begin{Definition}
	We say a diagonalisable endomorphism $K$ is {\it integrable}, if around any point in some open and dense set we f\/ind local
	coordinates such that the coordinate vectors are eigenvectors of $K$.
\end{Definition}

The above geometric def\/inition of integrability can be cast into a system of partial dif\/ferential equations, involving the
Nijenhuis torsion of $K$ def\/ined by
\begin{gather}
    \label{eq:Nijenhuis}
	N\indices{^\alpha_{\beta\gamma}}
	=K\indices{^\alpha_\delta}\nabla_{[\gamma}K\indices{^\delta_{\beta]}}
	+\nabla_\delta K\indices{^\alpha_{[\gamma}}K\indices{^\delta_{\beta]}}	 .
\end{gather}
The following result was proven by Nijenhuis for endomorphisms with simple eigenvalues \cite{Nijenhuis}, but is true in general~\cite{Schoebel}.

\begin{Proposition}
	A diagonalisable endomorphism field $K\indices{^\alpha_\beta}$ on a Riemannian manifold is integrable if and only if it
	satisfies the Nijenhuis integrability conditions
	\begin{subequations}
		\label{eq:TNS}
		\begin{gather}
			0=N\indices{^\delta_{[\beta\gamma}}g_{\alpha]\delta},\\
			0=N\indices{^\delta_{[\beta\gamma}}K_{\alpha]\delta},\\
			0=N\indices{^\delta_{[\beta\gamma}}K_{\alpha]\varepsilon}K\indices{^\varepsilon_\delta},
		\end{gather}
	\end{subequations}
	where the square brackets denote complete antisymmetrisation in the enclosed indices.
\end{Proposition}
At the end of the 19th century, Paul St\"ackel showed that any system of orthogonal separation coordinates gives rise to what we
now call a \emph{St\"ackel system}.  And in 1934, Luther P.~Eisenhart proved the converse.  This established a one-to-one
correspondence between orthogonal separation coordinates and St\"ackel systems in the following sense.
\begin{Definition}
	A {\it St\"ackel system} on an $n$-dimensional Riemannian manifold is an $n$-dimensional space of integrable Killing tensors
	which mutually commute in the algebraic sense\footnote{Eisenhart's original def\/inition of a St\"ackel system was more restrictive, but is equivalent to the one given here~\cite{Benenti93}.}\footnote{There is an alternative notion of ``commuting'' for Killing tensors, namely with respect to the Poisson bracket, which
		is dif\/ferent but closely related for a St\"ackel system.  Indeed, one can replace the condition that the Killing tensors
		in a St\"ackel system are integrable by the condition that they mutually commute under the Poisson bracket~\cite{Benenti93}.  In this article the term ``commuting'' will always be used in the algebraic sense.}.
\end{Definition}

\begin{Theorem}[\protect{\cite{Eisenhart, Staeckel91}}]
	On a Riemannian manifold, there is a bijective correspondence between orthogonal separation coordinates and St\"ackel systems.
\end{Theorem}

Benenti later proved that a St\"ackel system can be represented by a single integrable Killing tensor with simple eigenvalues.
More precisely:
\begin{Theorem}[\cite{Benenti93}]\quad
	\begin{enumerate}\itemsep=0pt
		\item[$1.$] Every St\"ackel system contains a Killing tensor with simple eigenvalues.
		\item[$2.$] Conversely, every integrable Killing tensor with simple eigenvalues is contained in some St\"ackel system.
		\item[$3.$] A St\"ackel system is uniquely determined by such a Killing tensor.
	\end{enumerate}
\end{Theorem}

In particular, the separation coordinates can be recovered from a generic representative in a St\"ackel system by solving the
eigenproblem and integrating the distributions normal to the eigenvectors.

Although initially introduced to classify separation coordinates, St\"ackel systems play a role far beyond separation of
variables, as they constitute an important class of completely integrable dynamical systems.  The construction and
classif\/ication of St\"ackel systems in constant curvature is, for instance, equivalent to the construction and classif\/ication of
Gaudin magnets \cite{Kalnins&Kuznetsov&Miller,Kuznetsov92b, Kuznetsov92a}.  St\"ackel systems also arise from Killing--Yano towers
and in black hole integrable models \cite{Frolov&Zelnikov}.  Furthermore, there exists a natural quantisation of Killing tensors
such that the quantum version of St\"ackel systems are the maximal commutative subalgebras of second order symmetries of the
Schr\"odinger operator~\cite{Benenti&Chanu&Rastelli}.

{\sloppy The Nijenhuis integrability conditions \eqref{eq:TNS} are homogeneous algebraic equations in~$K$ and~$\nabla K$, which
are invariant under isometries.  Separation of variables therefore naturally leads to the following problem.

}

\begin{Problem}
	\label{prob:solve}
	Determine the projective variety $\mathcal K(M)$ of all integrable Killing tensors on a~given Riemannian manifold~$M$.  That
	is, solve explicitly the Nijenhuis integrability conditions~\eqref{eq:TNS} for Killing tensors~\eqref{eq:Killing}.
	Moreover, find natural isometry invariants characterising the integrability of a Killing tensor.
\end{Problem}

Of course, there are many examples of separation coordinates and in constant curvature we even have a complete classif\/ication.
In principle this gives a description of the corresponding St\"ackel systems and integrable Killing tensors, although in practice
they are often obtained by intricate limiting processes.  Nevertheless, this only yields a description of $\mathcal K(M)$ as a
\emph{set}, whereas in Problem~\ref{prob:solve} we seek to elucidate the natural algebraic geometric structure of this set.
Note that in dimension two the Nijenhuis integrability conditions are void, so that $\mathcal K(M)$ is simply the projective
space of Killing tensors.  This means that Problem~\ref{prob:solve} is trivial in dimension two once the space of Killing
tensors is known.  However, to the best of our knowledge, already in dimensions three nothing is known about the set $\mathcal
K(M)$ as a projective variety, apart from trivial cases, not even for Euclidean space ${\mathbb R}^3$ or the sphere $S^3$ -- the two
examples for which Eisenhart already derived the complete list of separation coordinates.  Moreover, an explicit solution of the
equations \eqref{eq:TNS} has been considered intractable by various experts in this domain \cite{Horwood&McLenaghan&Smirnov05}.

Although there exists no general construction for St\"ackel systems, there is one which yields a large family of St\"ackel systems.
This family is interesting in itself, because it appears under dif\/ferent names in a number of dif\/ferent guises in the theory of
integrable systems:  as ``Newtonian systems of quasi-Lagrangian type''~\cite{Rauch-Wojciechowski&Marciniak&Lundmark}, as
``systems admitting special conformal Killing ten\-sors''~\cite{Crampin03a}, as ``cofactor systems'' \cite{Lundmark}, as
``bi-Hamiltonian structures''~\cite{Blaszak, Ibort&Magri&Marmo}, as ``bi-quasi-Hamiltonian systems''~\cite{Crampin03b, Crampin&Sarlet02}, as ``$L$-systems'' \cite{Benenti05} or as ``Benenti systems''~\cite{Bolsinov&Matveev}.  We
follow the latter nomenclature of Bolsinov and Matveev, because their description is best suited to our context.  This brings us
to the following question.
\begin{Question}	\label{Q:Benenti}
	Do Benenti systems parametrise the entire variety of integrable Killing tensors?
\end{Question}

Moreover, by construction, the set of Benenti systems is not only invariant under isometries, but even under the projective
group.  So, if the answer to Question~\ref{Q:Benenti} is ``no'' and the isometry group is a proper subgroup of the
projective group, the following question arises.
\begin{Question}
	\label{Q:invariance}
	Is the variety of integrable Killing tensors invariant under the projective group?
\end{Question}

Obviously, there are three groups acting on the variety $\mathcal K(M)$ of integrable Killing tensors:
\begin{enumerate}\itemsep=0pt
	\item 
${\mathbb R}$, acting by addition of multiples of the metric $g$.
	\item 
The isometry group $\operatorname{Isom}(M)$.
	\item 
${\mathbb R}^*$, acting by multiplication.
\end{enumerate}

The actions~1 
and~3 
result from the fact that adding the identity to an endomorphism or
multiplying it by a constant does not alter the eigenspaces and hence does not af\/fect integrability.  The actions~1 
and~2 
commute and the action~3 
descends to the quotient.  Hence the
essential information about integrable Killing tensors is encoded in the projectivisation of the quotient
\begin{gather}
	\label{eq:quotient}
	\mathcal K(M)/\bigl(\operatorname{Isom}(M)\times{\mathbb R} g\bigr)	 .
\end{gather}
This leads to the following problem.
\begin{Problem}
	\label{prob:quotient}
	Determine explicitly the quotient \eqref{eq:quotient}, which describes integrable Killing tensors modulo obvious
	symmetries.
\end{Problem}

Recall that orthogonal separation coordinates are in bijective correspondence with St\"ackel systems.  Thus, in view of the
classif\/ication problem for orthogonal separation coordinates, one is interested in linear subspaces of commuting integrable
Killing tensors, rather than in the integrable Killing tensors themselves.  This leads to the following problem.
\begin{Problem}
	\label{prob:Staeckel}
	Determine explicitly the space of all St\"ackel systems on a given Riemannian manifold $M$ $($modulo isometries$)$.  That is,
	find all maximal linear subspaces of commuting Killing tensors in $\mathcal K(M)$.
\end{Problem}

The separation coordinates can be recovered from an integrable Killing tensor with simple eigenvalues, because the coordinate
hypersurfaces are orthogonal to the eigenvectors.  This leads to the f\/inal step in the determination of separation coordinates
via integrable Killing tensors and St\"ackel systems.
\begin{Problem}
	\label{prob:coordinates}
	Find a $($canonical$)$ representative in each St\"ackel system for which it is possible to solve the eigenvalue problem and which
	has simple eigenvalues.  Deduce the classification of orthogonal separation coordinates and give natural isometry
	invariants.
\end{Problem}

In this article we will completely solve Problems~\ref{prob:solve}, \ref{prob:quotient}, \ref{prob:Staeckel} and \ref{prob:coordinates} for the three-dimensional sphere~$S^3$ equipped with the round metric
and give answers to Questions~\ref{Q:Benenti} and~\ref{Q:invariance}.  This is done without relying on any computer algebra.  One of the reasons we choose~$S^3$ is that it has the simplest isometry group among the $3$-dimensional constant curvature manifolds.  Our solution will
reveal diverse facets of the problem of separation of variables, related to representation theory, algebraic geometry,
geometric invariant theory, geodesic equivalence, Stashef\/f polytopes or moduli spaces of stable algebraic curves.

\subsection{Prior results}

In his fundamental paper Eisenhart also derived a list of orthogonal separation coordinates for~${\mathbb R}^3$ and~$S^3$~\cite{Eisenhart}.  Separation coordinates on $S^3$ were also classif\/ied in \cite{Olevskii}.  These results have been generalised
by Kalnins and Miller, who gave an iterative diagrammatical procedure for the construction of orthogonal separation coordinates
on spaces of constant curvature in any dimension~\cite{Kalnins,Kalnins&Miller86}.  Furthermore, Eisenhart computed the St\"ackel
systems on ${\mathbb R}^3$ in their respective separation coordinates.  For~$S^3$ the quantum point of view was adopted in~\cite{Kalnins&Miller&Winternitz}.

In dimension two virtually everything is known.  Koenigs classif\/ied all spaces with at least three independent Killing tensors.
Apart from constant curvature spaces these are four so-called Darboux spaces and eleven multiparameter families of Koenigs
spaces \cite{Koenigs}.  For Darboux spaces the separation coordinates have been classif\/ied in
\cite{KKMW, Kalnins&Kress&Winternitz} and for Koenigs spaces this is straightforward.  Darboux spaces have analogs in higher
dimensions~\cite{BEHR}, which can be obtained from constant curvature spaces via a St\"ackel transform, and for the three
dimensional version it is also straightforward to f\/ind all separation coordinates.  Finally, superintegrability theory has led
to the construction of non-conformally f\/lat systems which are separable~\cite{Kalnins&Kress&Miller}.

Despite these advances, the topology and the geometry of the space of integrable Killing tensors or its quotient by isometries
has never been studied for any non-trivial case.  Probably the only advance towards an algebraic geometric study of these spaces
has been the construction of isometry invariants for Killing tensors which discriminate the dif\/ferent orthogonal separation
coordinates in three dimensional Euclidean space \cite{Horwood&McLenaghan&Smirnov05}.  This approach had not become viable until
recently, because it relies on an extensive use of modern computer algebra.  The algorithm, which depends on an a priori
knowledge of the separation coordinates, involves the solution of a linear system of approximately 250,000 equations in 50,000
unknowns and the complexity of the algorithm renders it impractical in dimensions greater than three \cite{DHMcLS}.  Nothing
similar is known for the sphere or hyperbolic space.  We are not aware of any applications of the resulting invariants to
elucidate the structure of the resulting quotient.

\subsection{Method and results}

Our method is based on a purely algebraic description of the vector space of Killing tensors on constant curvature manifolds in
combination with the corresponding algebraic equivalent of the Nijenhuis integrability conditions.  We consider the standard
$3$-sphere, isometrically embedded as the unit sphere in a $4$-dimensional Euclidean vector space $(V,g)$:
\[
	S^3=\{x\in V\colon\, g(x,x)=1\}\subset V.
\]
The vector space of Killing tensors $K$ on $S^3$ is then naturally isomorphic to the vector space of algebraic curvature tensors on the
ambient space $V$ \cite{McLenaghan&Milson&Smirnov}.  These are (constant) tensors $R_{a_1b_1a_2b_2}$ on $V$ having the
symmetries of a Riemannian curvature tensor.  The isomorphism is explicitly given by
\begin{gather}
	\label{eq:correspondence}
	K_x(v,w)\coloneq R_{a_1b_1a_2b_2}x^{a_1}x^{a_2}v^{b_1}w^{b_2},
	\qquad
	x\in S^3,\qquad
	v,w\in T_xS^3	 ,
\end{gather}
and equivariant with respect to the natural actions of the isometry group $\OG(V,g)$ on both spaces.  Under this correspondence
the Nijenhuis integrability conditions~\eqref{eq:TNS} for a Killing tensor are equivalent to the following algebraic integrability
conditions for the associated algebraic curvature tensor~\cite{Schoebel}:
\begin{subequations}
	\label{eq:AIC}
	\begin{gather}
		\label{eq:AICI}
		\young(\atwo,\btwo,\ctwo,\dtwo)
		g_{ij}
		R\indices{^i_{b_1a_2b_2}}
		R\indices{^j_{d_1c_2d_2}}
		 =0,\\ \nonumber\\
		\label{eq:AICII}
		\smash[t]{\young(\atwo,\btwo,\ctwo,\dtwo)}
		\young(\aone \bone \cone \done)
		g_{ij}
		g_{kl}
		R\indices{^i_{b_1a_2b_2}}
		R\indices{^j_{a_1}^k_{c_1}}
		R\indices{^l_{d_1c_2d_2}}
		 =0		 .
	\end{gather}
\end{subequations}
Here the Young symmetrisers on the left hand side denote complete antisymmetrisation in the indices $a_2$, $b_2$, $c_2$, $d_2$
respectively complete symmetrisation in the indices $a_1$, $b_1$, $c_1$, $d_1$.  These homogeneous algebraic equations def\/ine the space
of integrable Killing tensors as a projective variety.

The key to the solution of the algebraic integrability conditions is a reinterpretion in terms of the decomposition of algebraic
curvature tensors into irreducible representations of the isometry group~-- the selfdual and anti-selfdual Weyl part, the
trace-free Ricci part and the scalar part.  Under the action of the isometry group these components can be brought to a certain
normal form, which simplif\/ies considerably if the algebraic integrability conditions are imposed.

We regard an algebraic curvature tensor on $V$ as a symmetric endomorphism on $\Lambda^2V$.  In this sense we will say that an
algebraic curvature tensor $R$ is {\it diagonal} if $R_{ijkl}=0$ unless $\{i,j\}=\{k,l\}$.  In dimension four a diagonal
algebraic curvature tensor can be parametrised as
\begin{gather}
	\label{eq:diagonal}
	R_{0\alpha0\alpha} =w_\alpha+t_\alpha+\tfrac s{12}, \qquad
	R_{\beta\gamma\beta\gamma} =w_\alpha-t_\alpha+\tfrac s{12}	 ,
\end{gather}
where $(\alpha,\beta,\gamma)$ denotes any cyclic permutation of $(1,2,3)$.  In this parametrisation the $w_\alpha$ are the
eigenvalues of the self-dual Weyl part, which turn out to be the same as those of the anti-self-dual Weyl part, the $t_\alpha$
parametrise the eigenvalues of the trace free Ricci tensor and $s$ is the scalar curvature.  In a f\/irst step we prove that
diagonal algebraic curvature tensors def\/ine a~(local) slice for the action of the isometry group on integrable Killing tensors.
\begin{Theorem}[integrability implies diagonalisability]
	\label{thm:slice}
	Under the action of the isometry group any integrable Killing tensor on $S^3$ is equivalent to one with a diagonal
	algebraic curvature tensor.
\end{Theorem}

The f\/irst algebraic integrability condition is identically satisf\/ied on diagonal algebraic curvature tensors and the second
reduces to a single algebraic equation on the diagonal entries.  This yields a description of integrable Killing tensors with
diagonal algebraic curvature tensor as a~projective variety, equipped with the residual action of the isometry group.  We call
this variety the {\it Killing--St\"ackel variety $($KS-variety$)$}.
\begin{Theorem}[integrability in the diagonal case]
	\label{thm:det}
	A Killing tensor on $S^3$ with diagonal algebraic curvature tensor \eqref{eq:diagonal} is integrable if and only if the
	associated matrix
	\begin{gather}
		\label{eq:matrix}
		\begin{pmatrix}
			\Delta_1&-t_3&t_2\\
			t_3&\Delta_2&-t_1\\
			-t_2&t_1&\Delta_3
		\end{pmatrix}
	\end{gather}
	has zero trace and zero determinant, where
		$\Delta_1 \coloneq w_2-w_3$,
		$\Delta_2 \coloneq w_3-w_1$,
		$\Delta_3 \coloneq w_1-w_2$.
	This defines a projective variety in $\mathbb P^4$ which carries a natural $S_4$-action, given by conjuga\-ting~\eqref{eq:matrix} with the symmetries of the regular octahedron in~$\mathbb R^3$ with vertices~$\pm e_i$ and adjacent faces
	oppositely oriented.  Two diagonal algebraic curvature tensors are equivalent under the isometry group if and only if their
	associated matrices~\eqref{eq:matrix} are equivalent under this $S_4$-action.
\end{Theorem}

Combining Theorems \ref{thm:slice} and \ref{thm:det} yields a solution to Problem~\ref{prob:quotient}:  The
projectivisation of the quotient \eqref{eq:quotient} is isomorphic to the quotient of the KS-variety by its natural
$S_4$-action.  Moreover, in the course of the proof of both theorems we establish isometry invariants for the integrability of a~Killing tensor.
\begin{Proposition}[isometry invariants]
	\label{prop:invariants}
	Let $\mathbf W$ and $\mathbf T$ be the Weyl and Ricci part of the algebraic curvature tensor corresponding to a Killing
	tensor on $S^3\subset V$, regarded as endomorphisms on $\Lambda^2V$.  Then this Killing tensor is integrable if and only if
	\begin{enumerate}\itemsep=0pt
		\item[$1)$] $\tr\,[\mathbf W,\mathbf T]^2=0$,
		\item[$2)$] $\tr (*\mathbf W^2)=\tr (*\mathbf W^3)=0$, where $\Hodge$ denotes the Hodge operator on $\Lambda^2V$ and
		\item[$3)$] the endomorphisms $\mathbf I$, $\mathbf W$ and $\mathbf W^2-\mathbf T^2$ are linearly dependent.
	\end{enumerate}
\end{Proposition}

With this result the check whether a given Killing tensor on $S^3$ is integrable becomes algorithmically trivial.

\looseness=-1
Na\"{\i}vely, one can f\/ind the St\"ackel systems by solving the commutator equation $[K,\tilde K]=0$ for a~f\/ixed integrable Killing
tensor $K$.  But a far simpler way is to exploit the algebraic geometry of the KS-variety.  Recall that St\"ackel systems are
linear spaces of dimension $n$ and always contain the metric.  We show that the Killing tensors in a St\"ackel system have
simultaneously diagona\-lisable algebraic curvature tensors.  Therefore St\"ackel systems on $S^3$ correspond to certain projective
lines in the KS-variety.  In general, the projective lines in a projective variety constitute themselves a projective variety,
called the {\it Fano variety}.  That is, we are looking for a subvariety of the Fano variety of the KS-variety.  By def\/inition,
the KS-variety is a linear section of the variety of $3\times3$ matrices with vanishing determinant, the so called full
determinantal variety.  Fortunately, the Fano variety of a full determinantal variety is well understood just in the case we
need here.

So we only have to use the commutator equation to \emph{check} which projective lines correspond to St\"ackel systems, rather than
to actually \emph{solve} it.  In order to do so, we derive an algebraic equivalent of the commutator equation in terms of the
corresponding algebraic curvature tensors.  This entails a detailed description of the algebraic geometry of the KS-variety and
the St\"ackel systems therein.  In particular, we can relate the algebraic geometric properties of the KS-variety to the geometric
properties of the corresponding integrable Killing tensors.  For this, we remark that there is a natural way how to extend
(integrable) Killing tensors from an embedded sphere $S^m\subset S^n$ to the ambient sphere $S^n$ for $m<n$.
\begin{Theorem}
	\label{thm:KS-variety}
	For $M=S^3$ the projectivisation of the quotient \eqref{eq:quotient} is isomorphic to the quotient of the KS-variety by its
	natural action of $S_4$.  Furthermore:
	\begin{enumerate}\itemsep=0pt
		\item[$1.$]
			The singularities of the KS-variety correspond to extensions of the metric on $S^2\subset S^3$ respectively on
			$S^1\subset S^3$.
		\item[$2.$]
			St\"ackel systems determine projective lines in the KS-variety which consist of matrices \eqref{eq:matrix}
			annihilating a given vector.
		\item[$3.$]
			Projective planes in the KS-variety that consist of matrices~\eqref{eq:matrix} annihilating a fixed vector
			correspond to extensions of $($integrable$)$ Killing tensors from~$S^2$ to~$S^3$.
		\item[$4.$]
			The projective plane of antisymmetric matrices in the KS-variety corresponds to the projective space of special
			Killing tensors.  These are Killing tensors $K$ of the form
			\[
				K=L-(\tr L)g,
			\]
			where $L$ is a special conformal Killing tensor\/\footnote{Cf.\ Def\/inition~\ref{def:SCKT}.}.
	\end{enumerate}
\end{Theorem}

Based on the algebraic description of St\"ackel systems we derive an algebraic description of Benenti systems as well and answer
Question~\ref{Q:Benenti}.  The answer is ``almost'', in the following sense.  On one hand, there are St\"ackel systems on
$S^3$, which are not Benenti systems.  So the answer is clearly ``no''.  But on the other hand these are exactly the extensions
of Killing tensors from $S^2$ and for $S^2$ the answer is easily seen to be ``yes''.

From the algebraic description of Benenti systems we derive the following parametrisation of~$\mathcal K(S^3)$, which is
equivariant under the isometry group.  This solves Problem~\ref{prob:solve}.
\begin{Theorem}[equivariant parametrisation]
	\label{thm:parametrisation}
	The algebraic variety $\mathcal K(S^3)$ is the closure of the set of algebraic curvature tensors of the form
	\[
		\lambda_2h\varowedge h+\lambda_0g\varowedge g,
		 \qquad
		h\in\Sym^2(V),
		\qquad
		\lambda_0,\lambda_2\in\mathbb R.
	\]
	Here $\varowedge$ denotes the Kulkarni--Nomizu product\/\footnote{Cf.~Def\/inition~\ref{def:Kulkarni-Nomizu}.} and
	$\frac12g\varowedge g$ is the algebraic curvature tensor of the metric.  That is, $\mathcal K(S^3)$ is the closure of the
	projective cone from the metric over the image of the Kulkarni--Nomizu square $h\mapsto h\varowedge h$, a map from symmetric
	tensors to algebraic curvature tensors.  Remark that the apex lies itself in this image.
\end{Theorem}
Since the set of algebraic curvature tensors of the form $h\varowedge h$ is invariant under the projective group $\mathbb
P\GLG(V)$ of $S^3$, the answer to Question~\ref{Q:invariance} is ``almost'' as well.  The reason is that although~${\mathcal
K}(S^3)$ is not projectively invariant, this ``invariant part'' still contains the essential part of each St\"ackel system unless
it is an extension from~$S^2$.  But again, on~$S^2$ the answer is easily seen to be ``yes''.  Our answers to Questions~\ref{Q:Benenti} and~\ref{Q:invariance}
indicate that one should think of the space of integrable Killing tensors on the sphere as composed of strata, where the
non-generic strata come from lower dimensional spheres.  Then each stratum is generated by Benenti systems and invariant under
the projective group in the respective dimension.

As a simple consequence of Theorem~\ref{thm:KS-variety} we get a canonical representative in each St\"ackel system.
\begin{Corollary}
	Each St\"ackel system on $S^3$ contains an essentially unique special Killing tensor.
\end{Corollary}
If the eigenvalues of the corresponding special conformal Killing tensor are simple, they def\/ine the respective separation
coordinates \cite{Crampin03a}.  Moreover, every special conformal Killing tensor $L$ on $S^n$ is the restriction of a constant
symmetric tensor $\hat L$ on the ambient space $V$ \cite{Matveev&Mounoud}.  This is the reason why we can solve the eigenvalue
problem for integrable Killing tensors explicitly.  In this manner we get a purely algebraic classif\/ication of separation
coordinates, which solves Problem~\ref{prob:coordinates}.  The details will be given later.

A solution to Problem~\ref{prob:Staeckel} can be derived from the description of St\"ackel systems in the KS-variety.  This
gives the classif\/ication space for separation coordinates on $S^3$ modulo isometries.
\begin{Theorem}[moduli space of separation coordinates]
	\label{thm:blow-up}
	The space of St\"ackel systems on~$S^3$ modulo isometries is homeomorphic to the quotient of the blow-up of~$\mathbb P^2$ in
	the four points $(\pm1\!:\!\pm1\!:\!\pm1)\in\mathbb P^2$ under the natural action of the symmetry group~$S_4$ of a
	tetrahedron in $\mathbb R^3$ inscribed in a cube with vertices $(\pm1,\pm1,\pm1)$.
\end{Theorem}
The fundamental domain of the $S_4$-action is a pentagon bounded by the non-generic separation coordinates.  In fact, the
multiplicities of the eigenvalues of the constant symmetric tensor~$\hat L$ on~$V$ identify this pentagon to the associahedron
(or Stashef\/f polytope) $K_4$ \cite{Stasheff63}.  These observations suggests a relation between separation coordinates on
spheres and moduli spaces of stable curves of genus zero with marked points \cite{Deligne&Mumford,Knudsen,Knudsen+}.  Indeed, our
thorough analysis of the case $S^3$ presented here has led to the following generalisation of Theorem~\ref{thm:blow-up}~\cite{Schoebel&Veselov}.
\begin{Theorem}	\label{thm:moduli}\sloppy
	The St\"ackel systems on $S^n$ with diagonal algebraic
	curvature tensor form a smooth projective variety isomorphic to
	the real Deligne--Mumford--Knudsen moduli space $\bar M_{0,n+2}(\mathbb R)$
	of stable genus zero curves with $n+2$ marked points.
\end{Theorem}

\subsection{Generalisations}

Most of our results are stated for constant curvature manifolds of arbitrary dimension.  What impedes a straightforward
generalisation to higher dimensional spheres is that in the proof of diagonalisability we make essential use of the Hodge
decomposition which only exists in dimension four.  Theorem~\ref{thm:moduli} together with a proof of the following
conjecture would show, that the classif\/ication space for separation coordinates on $S^n$ is the associahedron $K_{n+1}$.
\begin{Conjecture}\looseness=-1
	Under the action of the isometry group any integrable Killing tensor on $S^n$ is equivalent to one with a diagonal algebraic
	curvature tensor.
\end{Conjecture}
In Lorentzian signature the Hodge star squares to minus one.  A generalisation to hyperbolic space therefore naturally requires
complexif\/ication.  But over the complex numbers, all signatures are equivalent.   This indicates that one should solve the
algebraic integrability conditions over the complex numbers.  The variety of complex solutions then contains not only the
solution for hyperbolic space, but for all (non-f\/lat) pseudo-Riemannian constant curvature spaces as real subvarieties, singled
out by appropriate choices of a real structure.  In dimension three there is an elegant way to do this using spinors and the
Petrov classif\/ication of curvature tensors \cite{Milson&Schoebel}.  The main complication comes from the fact that
diagonalisability does not hold anymore in the complexif\/ied case.

A challenge remains without doubt a non-trivial algebraic geometric classif\/ication of integrable Killing tensors and orthogonal
separation coordinates on some non-constant curvature manifold.  This is out of scope of common, purely geometric techniques.
Note that the Nijenhuis integrability conditions \eqref{eq:TNS} can always be rewritten as homogeneous algebraic equations of
degree two, three and four on the space of Killing tensors, invariant under the isometry group action on Killing tensors.  This
endows the set of integrable Killing tensors on an arbitrary Riemannian manifold with the structure of a projective variety,
equipped with an action of the isometry group.  Our approach can therefore be extended to non-constant curvature as well.
Natural candidates are semi-simple Lie groups, homogeneous spaces or Einstein manifolds.  Notice that $S^3$ falls in either of
these categories.  Taking our solution for $S^3$ as a roadmap, we propose the following three-stage procedure:
\begin{enumerate}\itemsep=0pt
	\item
		Determine the space of Killing tensors and identify it explicitly as a representation of the isometry group (as in
		\cite{McLenaghan&Milson&Smirnov} for constant curvature manifolds).
	\item
		Translate the Nijenhuis integrability conditions for Killing tensors into purely algebraic integrability conditions (as
		in~\cite{Schoebel} for constant curvature manifolds).
	\item
		Solve the algebraic integrability conditions, study the algebraic geometric properties of the resulting variety and
		derive the complete classif\/ication of orthogonal separation coordinates (as in this article for~$S^3$).
\end{enumerate}
We consider the present work as a proof of concept that this algebraic geometric approach to the problem of separation of
variables is viable.

\subsection{Structure of the article}

The article is organised as follows:  In Section~\ref{sec:properties} we examine some properties of algebraic curvature
tensors that will be used to solve the algebraic integrability conditions in Section~\ref{sec:solution}.  We then describe
the algebraic geometry of this solution in Section~\ref{sec:algebraic} and identify the St\"ackel systems therein in
Section~\ref{sec:Staeckel}.  In Section~\ref{sec:geometric} we give various geometric constructions for integrable
Killing tensors and interpret them algebraically.  In the last section we demonstrate how the classif\/ication of separation
coordinates can be recovered within our framework.

\section{Properties of algebraic curvature tensors}
\label{sec:properties}

Algebraic curvature tensors on a vector space $V$ are four-fold covariant tensors $R_{a_1b_1a_2b_2}$ having the same algebraic
symmetries as a Riemannian curvature tensor:
\begin{subequations}
	\label{eq:R}
	\begin{alignat}{3}
		& \label{eq:R:anti}\text{antisymmetry:}\quad &&R_{b_1a_1a_2b_2} =-R_{a_1b_1a_2b_2}=R_{a_1b_1b_2a_2},& \\
		& \label{eq:R:pair}\text{pair symmetry:}\quad &&R_{a_2b_2a_1b_1} =R_{a_1b_1a_2b_2}, & \\
		& \label{eq:R:Bianchi}\text{Bianchi identity:}\quad &&R_{a_1b_1a_2b_2} +R_{a_1a_2b_2b_1}+R_{a_1b_2b_1a_2}=0.&
	\end{alignat}
\end{subequations}
Given a scalar product $g$ on $V$, we can raise and lower indices.  The symmet\-ries~\eqref{eq:R:anti} and~\eqref{eq:R:pair} then
allow us to regard an algebraic curvature tensor $R_{a_1a_2b_1b_2}$ on $V$ as a symmetric endomorphism
$R\indices{^{a_1b_1}_{a_2b_2}}$ on the space $\Lambda^2V$ of $2$-forms on $V$.  Since we will frequently change between both
interpretations, we denote endomorphisms by the same letter in boldface.

\subsection{Decomposition}

In the special case where $\dim V=4$, the Hodge star operator ``$\Hodge$'' def\/ines a decomposition
\[
	\Lambda^2V=\Lambda^2_+V\oplus\Lambda^2_-V
\]
of $\Lambda^2V$ into its $\pm1$ eigenspaces.  We can therefore write an algebraic curvature tensor $R$ and the Hodge star as
block matrices
\begin{subequations}
	\label{eq:branching}
	\begin{gather}
		\mathbf R
		=
		\left(
			\begin{array}{@{}c|c@{}}
				W_+&T_\mp\\\hline
				T_{\pm}&W_-
			\end{array}
		\right)
		+\frac s{12}
		\left(
			\begin{array}{@{}c|c@{}}
				\Id_+&0\\\hline
				0&\Id_-
			\end{array}
		\right), \qquad
		\Hodge =
		\left(
			\begin{array}{@{}c|c@{}}
				+\Id_+&0\\\hline
				0&-\Id_-
			\end{array}
		\right)
	\end{gather}
	with the $3\!\times\!3$-blocks satisfying
	\begin{gather}
		\label{eq:symmetric}
		W_+^\transpose=W_+,
		\qquad
		W_-^\transpose=W_-,
		\qquad
		T_\mp=T_\pm^\transpose,
		\\
		\label{eq:tracefree}
		\tr W_++\tr W_-=0,
		\\
		\label{eq:Bianchi}
		\tr W_+-\tr W_-=0 .
	\end{gather}
\end{subequations}
Here $\Id_+$ and $\Id_-$ denote the identity on $\Lambda^2_+V$ respectively $\Lambda^2_-V$.  The conditions \eqref{eq:symmetric}
assure symmetry, condition \eqref{eq:tracefree} says that $s=2\tr\mathbf R$ and condition \eqref{eq:Bianchi} is a reformulation
of the Bianchi identity \eqref{eq:R:Bianchi}.  Indeed, the Bianchi identity is equivalent to the vanishing of the
antisymmetrisation of $R_{a_1b_1a_2b_2}$ in all four indices.  In dimension four this can be written as
$
	\varepsilon^{a_2b_2a_1b_1}R_{a_1b_1a_2b_2}=0
$
or
\[
	\bm\varepsilon\indices{^{a_2b_2}_{a_1b_1}}\mathbf R\indices{^{a_1b_1}_{a_2b_2}}=0
	 ,
\]
where $\varepsilon^{a_2b_2a_1b_1}$ is the totally antisymmetric tensor.  Notice that
$\bm\varepsilon\indices{^{a_2b_2}_{a_1b_1}}$ is nothing else than the Hodge star operator and hence
\begin{gather}
	\label{eq:tr(*R)=0}
	\tr(\Hodge\mathbf R)=0	 ,
\end{gather}
where $\tr:\End(\Lambda^2V)\to{\mathbb R}$ is the usual trace.  Now \eqref{eq:Bianchi} is \eqref{eq:tr(*R)=0} applied to $\mathbf R$ in~\eqref{eq:branching}.

The space of algebraic curvature tensors is an irreducible $\GLG(V)$-representation and \eqref{eq:branching} gives a
decomposition of this representation into irreducible representations of the subgroup $\SOG(V,g)$ when $\dim V=4$.  As the
notation already suggests, we can relate these components to the familiar Ricci decomposition
\[
	\mathbf R=\mathbf W+\mathbf T+\mathbf S
\]
of an algebraic curvature tensor $R$ into
\begin{itemize}\itemsep=0pt
	\item a {\it scalar part}
		\[
			S_{a_1b_1a_2b_2}\coloneq\tfrac s{12}(g_{a_1a_2}g_{b_1b_2}-g_{a_1b_2}g_{b_1a_2})
		\]
		given by the {\it scalar curvature}
		\[
			s=g^{a_1a_2}g^{b_1b_2}R_{a_1b_1a_2b_2} ,
		\]
	\item a {\it trace free Ricci part}
		\[
			T_{a_1b_1a_2b_2}=\tfrac12(T_{a_1a_2}g_{b_1b_2}-T_{a_1b_2}g_{b_1a_2}-T_{b_1a_2}g_{a_1b_2}+T_{b_1b_2}g_{a_1a_2})
		\]
		given by the {\it trace free Ricci tensor}
		\[
			T_{a_1a_2}=g^{b_1b_2}R_{a_1b_1a_2b_2}-\tfrac s4g_{a_1a_2},
		\]
	\item
		and a {\it Weyl part}, given by the totally trace free {\it Weyl tensor} $W\coloneq R-T-S$.
\end{itemize}

It is not dif\/f\/icult to check that $[*,\mathbf T]=0$ and hence
\begin{gather}
	\label{eq:WTg}
	\mathbf W =
	\left(
		\begin{array}{@{}c|c@{}}
			W_+&0\\\hline
			0&W_-
		\end{array}
	\right), \qquad
	\mathbf T =
	\left(
		\begin{array}{@{}c|c@{}}
			0&T_\mp\\\hline
			T_{\pm}&0
		\end{array}
	\right),\qquad
	\mathbf S =
	\frac s{12}
	\left(
		\begin{array}{@{}c|c@{}}
			\Id_+&0\\\hline
			0&\Id_-
		\end{array}
	\right).
\end{gather}
We see that $W_+$ and $W_-$ are the self-dual and anti-self-dual part of the Weyl tensor.  Implicitly, the above
interpretation also provides an isomorphism
\begin{gather}
	\label{eq:iso}
	\begin{array}{@{}rcl}
		\Sym_0(V)&\stackrel{\isom}{\longrightarrow}&\operatorname{Hom}\big(\Lambda^2_+V,\Lambda^2_-V\big),\\
		T&\mapsto&T_\pm,
	\end{array}
\end{gather}
between trace-free Ricci tensors and homomorphisms from self-dual to anti-self-dual $2$-forms.

\subsection{The action of the isometry group}
\label{sub:groupaction}

The Killing tensor equation \eqref{eq:Killing} as well as the integrability conditions \eqref{eq:TNS} are invariant under the
action of the isometry group.  In other words, the variety of integrable Killing tensors is invariant under isometries.  This
will allow us later to put the algebraic curvature tensor of an integrable Killing tensor on $S^3$ to a certain normal form and
to solve the algebraic integrability conditions for this particular form.  Recall that the isometry group of $S^n\subset V$ is
$\OG(V)$ and that $\SOG(V)$ is the subgroup of orientation preserving isometries.

We f\/irst examine the induced action of $\SOG(V)$ on algebraic curvature tensors of the form~\eqref{eq:branching}.  The standard
action of $\SOG(V)$ on~$V$ induces a natural action on~$\Lambda^2V$ which is the adjoint action of~$\SOG(V)$ on its Lie algebra
under the isomorphism
\[
	\mathfrak{so}(V)\isom\Lambda^2V.
\]
Now consider the following commutative diagram of Lie group morphisms for $\dim V=4$:
\[
	\begin{CD}
		@.\SpinG(V)@>\isom>>\SpinG\big(\Lambda^2_+V\big)@.\times@.\SpinG\big(\Lambda^2_-V\big)\\
		@.@VV2:1V@VV2:1V@.@VV2:1V\\
		 \pi\colon@.\SOG(V)@>2:1>>\SOG\big(\Lambda^2_+V\big)@.\times@.\SOG\big(\Lambda^2_-V\big)@.\hookrightarrow\mspace{8mu}\SOG\big(\Lambda^2V\big)\mspace{24mu}\\
		@.U\makebox[0pt]{\hspace{68pt}$\longmapsto$}@.\mspace{64mu}(U_+@.,@.U_-)\mspace{64mu}@.\mapsto\begin{pmatrix}U_+&0\\0&U_-\end{pmatrix}
		 .
	\end{CD}
\]
Here the double covering $\pi$ in the second row is induced from the exceptional isomorphism in the f\/irst row via the universal
covering maps (vertical).  Under the induced isomorphism $\pi_*$ of Lie algebras,
\[
	\mathfrak{so}(V)\isom\mathfrak{so}\big(\Lambda^2_+V\big)\oplus\mathfrak{so}\big(\Lambda^2_-V\big),
\]
the adjoint action of an element $U\in\SOG(V)$ corresponds to the adjoint actions of $U_+\in\SOG(\Lambda^2_+V)$ and
$U_-\in\SOG(\Lambda^2_-V)$.  Since $\Lambda^2_+V$ has dimension three, the Hodge star operator gives an isomorphism
\[
	\mathfrak{so}\big(\Lambda^2_+V\big)\isom\Lambda^2\Lambda^2_+V\isom\Lambda^2_+V
\]
under which the adjoint action of $U_+$ corresponds to the standard action on $\Lambda^2_+V$ and similarly for $U_-$.  The Hodge
decomposition completes the above isomorphisms to a commutative diagram
\[
	\begin{CD}
		\pi_*\colon@.\mathfrak{so}(V)@>\isom>>\mathfrak{so}\big(\Lambda^2_+V\big)@.\;\oplus\;@.\mathfrak{so}\big(\Lambda^2_-V\big)\\
		@.@VV\isom V@VV\isom V@.@VV\isom V\\
		@.\Lambda^2V@>\isom>>\Lambda^2_+V@.\oplus@.\Lambda^2_-V
		 ,
	\end{CD}
\]
which shows that the natural action of $U\in\SOG(V)$ on $\Lambda^2V$ is given by
\begin{gather}
	\label{eq:blockdiagonal}
	\begin{pmatrix}
		U_+&0\\
		0&U_-
	\end{pmatrix}
	\in\pi\bigl(\SOG(V)\bigr)
	\subset\SOG\big(\Lambda^2V\big) .
\end{gather}
Hence the natural action of $U$ on $\End(\Lambda^2V)$ is given by conjugation with this matrix.  Restricting to algebraic
curvature tensors we get:
\begin{Proposition}[action of the isometry group]
	An element $U\in\SOG(V)$ acts on algebraic curvature tensors in the form \eqref{eq:branching} by conjugation with
	\eqref{eq:blockdiagonal}, i.e.\ via
	\begin{gather}
		\label{eq:trafos}
		W_+ \mapsto U_+^\transpose W_+U_+,\qquad
		W_-  \mapsto U_-^\transpose W_-U_-,\qquad
		T_\pm \mapsto U_-^\transpose T_\pm U_+
		 ,
	\end{gather}
	where $\pi(U)=(U_+,U_-)$ is the image of $U$ under the double covering
	\[
		\pi\colon \ \SOG(V)\to\SOG\big(\Lambda^2_+V\big)\times\SOG\big(\Lambda^2_-V\big)
		 .
	\]
\end{Proposition}
We can describe the twofold cover $\pi$ explicitly in terms of orthonormal bases on the spaces~$V$, $\Lambda^2_+V$ and
$\Lambda^2_-V$.  It maps an orthonormal basis $(e_0,e_1,e_2,e_3)$ of $V$ to the orthonormal bases
$(\eta_{+1},\eta_{+2},\eta_{+3})$ of $\Lambda^2_+V$ and $(\eta_{-1},\eta_{-2},\eta_{-3})$ of $\Lambda^2_-V$, def\/ined by
\begin{gather}
	\label{eq:Hodgebasis}
	\eta_{\pm\alpha}\coloneq\tfrac1{\sqrt2}(e_0\wedge e_\alpha\pm e_\beta\wedge e_\gamma)
\end{gather}
for each cyclic permutation $(\alpha,\beta,\gamma)$ of $(1,2,3)$.  From this description also follows that $\ker\pi=\{\pm\Id\}$.

The action of $\OG(V)$ on algebraic curvature tensors is now determined by the action of some orientation reversing element in
$\OG(V)$, say the one given by reversing the sign of $e_0$ and preserving $e_1$, $e_2$ and $e_3$.  This element maps
$\eta_{\pm\alpha}$ to $-\eta_{\mp\alpha}$.  Hence its action on algebraic curvature tensors in the form \eqref{eq:branching} is
given by conjugation with
\[
	\begin{pmatrix}
		0&-I\\
		-I&0
	\end{pmatrix}
\]
respectively by mapping
\begin{gather}
	\label{eq:orientationreversing}
	W_+ \mapsto W_-,\qquad
	T_\pm \mapsto T_\mp	 .
\end{gather}

\subsection{Aligned algebraic curvature tensors}
\label{sub:aligned}

$W_+$ and $W_-$ are symmetric and hence simultaneously diagonalisable under the action~\eqref{eq:trafos} of~$\SOG(V)$.  On the
other hand, the singular value decomposition shows that~$T_\pm$ is also diagona\-lisable under this action, although in general
not simultaneously with~$W_+$ and~$W_-$.  The following lemma gives a criterion when this is the case.
\begin{Lemma}
	Let $W_+$ and $W_-$ be symmetric endomorphisms on two arbitrary Euclidean vector spaces $\Lambda_+$ respectively $\Lambda_-$
	and suppose the linear map $T_\pm\colon\Lambda_+\to\Lambda_-$ satisf\/ies
	\[
		T_\pm W_+=W_-T_\pm
		 .
	\]
	Then there exist orthonormal bases for $\Lambda_+$ and $\Lambda_-$ such that $W_+$, $W_-$ and $T_\pm$ are
	simultaneously diagonal with respect to these bases.
\end{Lemma}
\begin{proof}
	It suf\/f\/ices to show that we can chose diagonal bases for $W_+$ and $W_-$ such that the matrix of $T_\pm$ has at most one
	non-zero element in each row and in each column, for the desired result can then be obtained by an appropriate permutation
	of the basis elements.  The above condition implies that $T_\pm$ maps eigenspaces of $W_+$ to eigenspaces of $W_-$ with the
	same eigenvalue.  Without loss of generality we can thus assume that $\Lambda_+$ and $\Lambda_-$ are eigenspaces of $W_+$
	respectively $W_-$ with the same eigenvalue.  But then $W_+$ and $W_-$ are each proportional to the identity and therefore
	diagonal in any basis.  In this case the lemma follows from the singular value decomposition of~$T_\pm$.
\end{proof}

\begin{Lemma}
	\label{lem:aligned}
	The following conditions are equivalent for an algebraic curvature tensor \eqref{eq:branching}:
\begin{enumerate} \itemsep=0pt
		\item[$1)$] 
$T_\pm W_+=W_-T_\pm$,
		\item[$2)$] 
$[\mathbf W,\mathbf T]=0$,
		\item[$3)$] 
$\mathbf W$ and $\mathbf T$ are simultaneously diagonalisable under~$\SOG(\Lambda^2V)$.
	\end{enumerate}
\end{Lemma}

\begin{proof}
	The equivalence of 1 and 2 follows from~\eqref{eq:WTg}.  For the equivalence of~2 and~3 it suf\/f\/ices to note that $\mathbf W$ and $\mathbf T$ are both symmetric and hence
	diagonalisable.
\end{proof}

The following terminology is borrowed from general relativity.
\begin{Definition}
	We say that an algebraic curvature tensor is {\it aligned}, if it satisf\/ies one of the equivalent conditions in
	Lemma~\ref{lem:aligned}.
\end{Definition}
By Lemma~\ref{lem:aligned}, for any aligned algebraic curvature tensor we can f\/ind an orthonormal basis of $V$ such that
\begin{subequations}
	\label{eq:aligned}
	\begin{gather}
		\label{eq:aligned:diagonal}
		W_+ =
		\left(
			\begin{matrix}
				w_{+1}&0&0\\
				0&w_{+2}&0\\
				0&0&w_{+3}
			\end{matrix}
		\right), \quad
		 W_- =
		\left(
			\begin{matrix}
				w_{-1}&0&0\\
				0&w_{-2}&0\\
				0&0&w_{-3}
			\end{matrix}
		\right),\quad
		 T_\pm =
		\left(
			\begin{matrix}
				t_1&0&0\\
				0&t_2&0\\
				0&0&t_3
			\end{matrix}
		\right)\!\!\!
	\end{gather}
	with
	\begin{gather}
		\label{eq:aligned:plane}
		w_{+1}+w_{+2}+w_{+3} =0,\qquad
		w_{-1}+w_{-2}+w_{-3} =0.
	\end{gather}
\end{subequations}
To simplify notation we agree that henceforth the indices $(\alpha,\beta,\gamma)$ will stand for an arbitrary cyclic permutation
of $(1,2,3)$.  Changing the basis in $\Lambda^2V$ from \eqref{eq:Hodgebasis} to $e_i\wedge e_j$, $0\le i<j\le4$, we obtain the
independent components of an aligned algebraic curvature tensor:
\begin{gather}
		R_{0\alpha0\alpha} =\frac{w_{+\alpha}+w_{-\alpha}}2+t_\alpha+\frac s{12},\qquad
		R_{\beta\gamma\beta\gamma} =\frac{w_{+\alpha}+w_{-\alpha}}2-t_\alpha+\frac s{12},\nonumber\\
		R_{0\alpha\beta\gamma} =\frac{w_{+\alpha}-w_{-\alpha}}2 .\label{eq:components}
	\end{gather}
The trace free Ricci tensor of this aligned algebraic curvature tensor is diagonal and given by
\begin{gather}
	\label{eq:iso:diagonal}
		T_{00} = t_\alpha+t_\beta+t_\gamma,\qquad 		T_{\alpha\alpha} = t_\alpha-t_\beta-t_\gamma,
\qquad
	t_\alpha =\frac{T_{00}+T_{\alpha\alpha}}2	 .
\end{gather}
This is nothing else than the restriction of the isomorphism \eqref{eq:iso} to diagonal tensors.

\subsection{Diagonal algebraic curvature tensors}

\begin{Definition}
	We call an algebraic curvature tensor on $V$ {\it diagonal} in a basis $e_i$ if it is diagonal as an element of
	$\End(\Lambda^2V)$ with respect to the associated basis $e_i\wedge e_j$, $i<j$.  We call it {\it diagonalisable} if it is
	diagonalisable under the adjoint action of $\SOG(V)$ on $\Lambda^2V$, i.e.\ under the subgroup
	$\pi\bigl(\SOG(V)\bigr)\subset\SOG(\Lambda^2V)$.
\end{Definition}
Of course, being a symmetric endomorphism on $\Lambda^2V$, an algebraic curvature tensor is always diagonalisable under the full
group $\SOG(\Lambda^2V)$.

\begin{Proposition}[diagonalisability criterion]
	\label{prop:diagonalisability}
	An algebraic curvature tensor is diagonalisable if and only if it is aligned and $W_+$ has the same characteristic
	polynomial as $W_-$.
\end{Proposition}
\begin{proof}
	Consider an aligned algebraic curvature tensor and suppose that $W_+$ has the same characteristic polynomial as $W_-$.  Then
	we can f\/ind a transformation \eqref{eq:trafos} such that, as $3\times3$ matrices, $W_+=W_-$.  The condition $W_-T_\mp=T_\pm
	W_+$ then implies that $T_\pm$ can be simultaneously diagonalised with $W_+=W_-$.  This means we can assume without loss of
	generality that
	\[
		w_{+\alpha}=w_{-\alpha}\eqcolon w_\alpha
	\]
	for $\alpha=1,2,3$ in \eqref{eq:aligned}, so that \eqref{eq:components} reads
	\begin{subequations}
		\label{eq:chiral}
		\begin{gather}
			R_{0\alpha0\alpha} =w_\alpha+t_\alpha+\frac s{12},\qquad
			R_{0\alpha\beta\gamma} =0 ,\qquad
			R_{\beta\gamma\beta\gamma} =w_\alpha-t_\alpha+\frac s{12}
		\end{gather}
		with $w_\alpha$ subject to
		\begin{gather}
			\label{eq:chiral:trace}
			w_1+w_2+w_3=0
			 .
		\end{gather}
	\end{subequations}
	Obviously this algebraic curvature tensor is diagonal.  On the other hand, any diagonal algebraic curvature tensor can be
	written in the form \eqref{eq:chiral}, because
		\begin{gather*}
			s =2 (R_{0101}+R_{2323}+R_{0202}+R_{3131}+R_{0303}+R_{1212}),\\
			t_\alpha =\frac{R_{0\alpha0\alpha}-R_{\beta\gamma\beta\gamma}}2,\qquad
			w_\alpha =\frac{R_{0\alpha0\alpha}+R_{\beta\gamma\beta\gamma}}2-\frac s{12}	.\tag*{\qed}
		\end{gather*} \renewcommand{\qed}{}
\end{proof}

\subsection{The residual action of the isometry group}
\label{sub:remainder}

For later use we need the stabiliser of the space of diagonal algebraic curvature tensors under the action of the isometry
group.  Such tensors are of the form \eqref{eq:aligned} with $w_{-\alpha}=w_{+\alpha}$ and are invariant under the
transformation \eqref{eq:orientationreversing}.  Hence it is suf\/f\/icient to consider only the subgroup $\SOG(V)$ of orientation
preserving isometries.

Let $U\in\SOG(V)$ with $\pi(U)=(U_+,U_-)$ be an element preserving the space of diagonal algebraic curvature tensors under the
action \eqref{eq:trafos}.  Then $U_+\in\SOG(\Lambda^2_+V)$ and $U_-\in\SOG(\Lambda^2_-V)$ preserve the space of diagonal
matrices on $\Lambda^2_+V$ respectively $\Lambda^2_-V$.  In the orthogonal group, the stabiliser of the space of diagonal
matrices under conjugation is the subgroup of signed permutation matrices.  They act by permuting the diagonal elements in
disregard of the signs.  In particular, in $\SOG(3)$ the stabiliser subgroup of the space of diagonal matrices on ${\mathbb R}^3$ is
isomorphic to $S_4$.  Under the isomorphism $S_4\isom S_3\ltimes V_4$, the permutation group $S_3$ is the subgroup of (unsigned)
permutation matrices and acts by permuting the diagonal elements, whereas the Klein four group $V_4$ is the subgroup of diagonal
matrices in $\SOG(3)$ and acts trivially.

Using the above, it is not dif\/f\/icult to see that the stabiliser of the space of diagonal algebraic curvature tensors under the
action \eqref{eq:trafos} is isomorphic to $S_4\times\ker\pi$, with $\ker\pi=\{\pm\Id\}$ acting trivially.  Under the isomorphism
$S_4\isom S_3\ltimes V_4$, the factor $S_3$ acts by simultaneous permutations of $(w_1,w_2,w_3)$ and $(t_1,t_2,t_3)$ and the
factor $V_4$ by simultaneous f\/lips of two signs in $(t_1,t_2,t_3)$.

\section{Solution of the algebraic integrability conditions}
\label{sec:solution}

This section is dedicated to proving Theorems~\ref{thm:slice} and \ref{thm:det}.

\subsection{Reformulation of the f\/irst integrability condition}
\label{sub:reformulation}

In the same way as for the Bianchi identity we can reformulate the f\/irst integrability condi\-tion~\eqref{eq:AICI} as
\[
	\mathbf R\indices{^{ib_1}_{a_2b_2}}
	\bm\varepsilon\indices{^{a_2b_2}_{c_2d_2}}
	\mathbf R\indices{^{c_2d_2}_{id_1}}
	=0
\]
or
\begin{gather}
	\label{eq:str(R*R)=0}
	r(\mathbf R\Hodge\mathbf R)=0,
\end{gather}
where
\begin{gather}
	\label{eq:str}
	\begin{array}{@{}rrcl}
		r\colon&\End(\Lambda^2V)&\longrightarrow&\End(V),\\
		&\mathbf E\indices{^{ij}_{kl}}&\mapsto&\mathbf E\indices{^{ij}_{kj}}
	\end{array}
\end{gather}
denotes the Ricci contraction.  Notice that $\mathbf E\coloneq\mathbf R\Hodge\mathbf R$ is symmetric, but does in general not
satisfy the Bianchi identity.  A proof of the following lemma can be found in \cite{Singer&Thorpe}.
\begin{Lemma}
	The kernel of the Ricci contraction \eqref{eq:str} is composed of endomorphisms $\mathbf E$ satisfying
	\begin{gather*}
		\Hodge\mathbf E\Hodge =\mathbf E^\transpose,\qquad
		\tr\mathbf E =0		 .
	\end{gather*}
\end{Lemma}
Applying this to the symmetric endomorphism
\[
	\mathbf E
	\coloneq\mathbf R\Hodge\mathbf R
	=\left(
		\begin{array}{@{}cc@{}}
			\multirow{2}{*}{$(W_++\tfrac s{12}\Id_+)^2-T_\mp T_{\pm}$}
			&\multicolumn{1}{|l}{(W_++\tfrac s{12}\Id_+)T_\mp}\\
			&\multicolumn{1}{|r}{-T_\mp(W_-+\frac s{12}\Id_-)}
			\\\hline
			\multicolumn{1}{l|}{T_{\pm}(W_++\tfrac s{12}\Id_+)}&\multirow{2}{*}{$T_\pm T_\mp-(W_-+\tfrac s{12}\Id_-)^2$}\\
			\multicolumn{1}{r|}{-(W_-+\frac s{12}\Id_-)T_\pm}
		\end{array}
	\right)
	 ,
\]
in \eqref{eq:str(R*R)=0} and using \eqref{eq:tracefree} we get
\begin{subequations}
	\label{eq:AICI:dec}
	\begin{gather}
		W_+T_\mp =T_\mp W_-,\\
		\label{eq:squaretrace}
		\tr W_+^2 =\tr W_-^2.
	\end{gather}
\end{subequations}
This shows:
\begin{Proposition}
	\label{prop:AICI}
	The first algebraic integrability condition for an algebraic curvature tensor $\mathbf R=\mathbf W+\mathbf T$ with Weyl part
	$\mathbf W$ and $($not necessarily trace free$)$ Ricci part $\mathbf T$ is equivalent to
	\begin{gather*}
		[ \mathbf W,\mathbf T ] =0,\qquad
		\tr\bigl(\Hodge\mathbf W^2\bigr) =0 .
	\end{gather*}
	In particular, an algebraic curvature tensor satisfying the first algebraic integrability condition is aligned.
\end{Proposition}
As a consequence of Proposition~\ref{prop:diagonalisability} we get:
\begin{Corollary}
	\label{cor:diagonal->AICI}
	A diagonalisable algebraic curvature tensor satisfies the first algebraic integrability condition.
\end{Corollary}

\subsection{Integrability implies diagonalisability}
\label{sub:AICI}

The aim of this subsection is to prove Theorem~\ref{thm:slice}.  In view of Propositions \ref{prop:diagonalisability} and
\ref{prop:AICI} it is suf\/f\/icient to show that if an algebraic curvature tensor satisf\/ies the algebraic integrability conditions,
then $W_+$ and $W_-$ have the same characteristic polynomial.  We will f\/irst prove that the f\/irst algebraic integrability condition
implies that $W_+$ has the same characteristic polynomial either as $+W_-$ or as $-W_-$.  We then prove that the latter
contradicts the second algebraic integrability condition.

We saw that an algebraic curvature tensor which satisf\/ies the f\/irst algebraic integrability condition is aligned and thus of the
form \eqref{eq:aligned} in a suitable orthogonal basis.  The f\/irst integrability condition in the form \eqref{eq:AICI:dec} then
translates to
\begin{subequations}
	\label{eq:123}
	\begin{gather}
		\label{eq:123:sphere}
		w_{+1}^2+w_{+2}^2+w_{+3}^2 =r^2,\qquad
		w_{-1}^2+w_{-2}^2+w_{-3}^2 =r^2
	\end{gather}
	for some $r\ge 0$ and
	\begin{gather}
		\label{eq:123:commuting}
		w_{+\alpha}t_\alpha=t_\alpha w_{-\alpha}.
	\end{gather}
\end{subequations}
If we regard $(w_{+1},w_{+2},w_{+3})$ and $(w_{-1},w_{-2},w_{-3})$ as vectors in ${\mathbb R}^3$, then each equation in
\eqref{eq:aligned:plane} describes the plane through the origin with normal $(1,1,1)$ and each equation in~\eqref{eq:123:sphere}
the sphere of radius~$r$ centered at the origin.  Hence the solutions to both equations lie on a circle and can be parametrised
in polar coordinates by a radius $r\ge0$ and angles $\varphi_+$, $\varphi_-$ as
\begin{gather}
	\label{eq:dreibeine}
	w_{+\alpha} =r\cos\bigl(\varphi_++\alpha\tfrac{2\pi}3\bigr),\qquad
	w_{-\alpha} =r\cos\bigl(\varphi_-+\alpha\tfrac{2\pi}3\bigr)
	 .
\end{gather}
Assume now that $W_+$ and $W_-$ have dif\/ferent eigenvalues.  Then we have $T_\pm=0$.  Indeed, if $t_\alpha\not=0$ for some
$\alpha$, then \eqref{eq:123:commuting} shows that $w_{+\alpha}=w_{-\alpha}$ for this index $\alpha$.  But then
\eqref{eq:dreibeine} implies that $W_+$ and $W_-$ have the same eigenvalues, which contradicts our assumption.  Hence $T_\pm=0$.

Without loss of generality we will set $s=0$ in what follows, because integrability is independent of the scalar curvature $s$.
This follows from the fact that the Killing tensor corresponding to the scalar component in the Ricci decomposition is
proportional to the metric and thus to the identity when regarded as an endomorphism.  But adding a multiple of the identity
only alters the eigenvalues, not the eigendirections, and therefore has no ef\/fect on integrability.

After these considerations, the independent components \eqref{eq:components} of our algebraic curvature tensor are
\begin{gather}
	\label{eq:achiral}
	R_{0\alpha0\alpha} =\frac{w_{+\alpha}+w_{-\alpha}}2=R_{\beta\gamma\beta\gamma},\qquad
	R_{0\alpha\beta\gamma} =\frac{w_{+\alpha}-w_{-\alpha}}2
	 .
\end{gather}

Choosing $a_1=b_1=c_1=d_1=0$ in the second integrability condition \eqref{eq:AICII} for this algebraic curvature tensor
yields
\[
	g^{ij}g^{kl}
	\young(\atwo,\btwo,\ctwo,\dtwo)
	R_{i0a_2b_2}
	R_{j0k0}
	R_{l0c_2d_2}
	=
	\sum_{\alpha=1,2,3}
	\young(\atwo,\btwo,\ctwo,\dtwo)
	R_{\alpha0a_2b_2}
	R_{\alpha0\alpha0}
	R_{\alpha0c_2d_2}
	=0
	 ,
\]
since only terms with $i=j=k=l\not=0$ contribute to the contraction.  For $(a_2,b_2,c_2,d_2)=(0,1,2,3)$ this is gives
\[
	\sum R_{\alpha00\alpha}R_{\alpha0\alpha0}R_{\alpha0\beta\gamma}
	=\sum(R_{0\alpha0\alpha})^2R_{0\alpha\beta\gamma}
	=0
	 ,
\]
where the sums run over the three cyclic permutations $(\alpha,\beta,\gamma)$ of $(1,2,3)$.  We conclude that
\[
	\sum_{\alpha=1,2,3}
	(w_{+\alpha}+w_{-\alpha})^2
	(w_{+\alpha}-w_{-\alpha})
	=0
	 .
\]
Substituting \eqref{eq:dreibeine} into this equation yields, after some trigonometry,
\[
	r^3\cos^2\frac{\varphi_+-\varphi_-}2\sin\frac{\varphi_+-\varphi_-}2\sin\frac32(\varphi_++\varphi_-)=0
	 .
\]
Substituting the solutions
\begin{gather*}
	r =0,\qquad
	\varphi_+ =\varphi_-+k\pi,\qquad
	\varphi_+ =-\varphi_-+k\tfrac{2\pi}3,\qquad
	k \in\mathbb Z
\end{gather*}
of this equation back into \eqref{eq:dreibeine} shows that the set $\{w_{-1},w_{-2},w_{-3}\}$ is either equal to
$\{w_{+1},w_{+2}$, $w_{+3}\}$ or to $\{-w_{+1},-w_{+2},-w_{+3}\}$.  The f\/irst case is excluded by assumption and we conclude that
$W_+$ and $-W_-$ have the same eigenvalues.  Consequently we can choose an orthonormal basis of $V$ such that
\begin{gather}
	\label{eq:opposed}
	w_{+\alpha}=-w_{-\alpha}\eqcolon w_\alpha
	 .
\end{gather}
Then the only remaining independent components of the algebraic curvature tensor \eqref{eq:achiral} are
\begin{gather}
	\label{eq:R0123}
	R_{0\alpha\beta\gamma}=w_\alpha
	 .
\end{gather}
Now consider the second integrability condition \eqref{eq:AICII} for $(a_2,b_2,c_2,d_2)=(0,1,2,3)$, written as
\begin{gather}
	\label{eq:tobecontracted}
	\young(\aone\bone\cone\done)
	g^{ij}g^{kl}E_{ib_1ld_1}R_{ja_1kc_1}=0
\end{gather}
with
\begin{gather}
	\label{eq:E}
	E_{ib_1ld_1}\coloneq
	\frac14\young(\atwo,\btwo,\ctwo,\dtwo)R_{ib_1a_2b_2}R_{ld_1c_2d_2}.
\end{gather}
The product $R_{ib_1a_2b_2}R_{ld_1c_2d_2}$ is zero unless $\{i,b_1,a_2,b_2\}=\{l,d_1,c_2,d_2\}=\{0,1,2,3\}$.  Since
$\{a_2,b_2,c_2,d_2\}=\{0,1,2,3\}$, we have $E_{ib_1ld_1}=0$ unless $\{i,b_1,l,d_1\}=\{0,1,2,3\}$.  The symmetries
\begin{gather}
	\label{eq:symmetries}
	E_{b_1ild_1} =-E_{ib_1ld_1}=E_{ib_1d_1l},\qquad
	E_{ib_1ld_1} =E_{ld_1ib_1}
\end{gather}
then imply that the only independent components of \eqref{eq:E} are
\begin{gather}
	\label{eq:E0123}
	E_{0\alpha\beta\gamma}=w_\alpha^2
	 .
\end{gather}
Now consider the full contraction of \eqref{eq:tobecontracted}:
\[
	g^{a_1b_1}g^{c_1d_1}
	\young(\aone\bone\cone\done)
	E\indices{^i_{b_1}^k_{d_1}}R_{ia_1kc_1}	=0	 .
\]
Since the Ricci tensor of \eqref{eq:R0123} is zero, this is results in
\[
	E^{ijkl}R_{ijkl}+E^{ijkl}R_{ilkj}=0	 .
\]
Using the symmetries \eqref{eq:symmetries} of $E_{ijkl}$ and those of $R_{ijkl}$, we can always permute a particular index to
the f\/irst position:
\[
	E^{0jkl}R_{0jkl}+E^{0jkl}R_{0lkj}=0
	 .
\]
And since both $E^{ijkl}$ and $R_{ijkl}$ vanish unless $\{i,j,k,l\}=\{0,1,2,3\}$, we get
\[
	\sum\big(E^{0\alpha\beta\gamma}R_{0\alpha\beta\gamma}+E^{0\gamma\beta\alpha}R_{0\gamma\beta\alpha}\big)+
	\sum\big(E^{0\alpha\beta\gamma}R_{0\gamma\beta\alpha}+E^{0\gamma\beta\alpha}R_{0\alpha\beta\gamma}\big)=0.
\]
Here again, the sums run over the three cyclic permutations $(\alpha,\beta,\gamma)$ of $(1,2,3)$.  Substitu\-ting~\eqref{eq:R0123}
and \eqref{eq:E0123} yields
\[
	\sum\big(w_\alpha^3+w_\gamma^3\big)-\sum\big(w_\alpha^2w_\gamma+w_\gamma^2w_\alpha\big)=0.
\]
With \eqref{eq:chiral:trace} the second term can be transformed to
\[
	\sum\big(w_\alpha^2w_\gamma+w_\gamma^2w_\alpha\big)=
	\sum\big(w_\beta^2w_\alpha+w_\beta^2w_\gamma\big)=
	\sum w_\beta^2(w_\alpha+w_\gamma)=
	-\sum w_\beta^3,
\]
implying
\[
	w_1^3+w_2^3+w_3^3=0.
\]
By \eqref{eq:chiral:trace} and Newton's identity
\begin{gather*}
	w_1^3+w_2^3+w_3^3
	=3w_1w_2w_3
	-3(w_1w_2+w_2w_3+w_3w_1)(w_1+w_2+w_3)
	+(w_1+w_2+w_3)^3
\end{gather*}
this is equivalent to
\[
	w_1w_2w_3=0.
\]
Using \eqref{eq:chiral:trace} once more, this shows that $\{w_1,w_2,w_3\}=\{-w,0,+w\}$ for some $w\in{\mathbb R}$.  But then
\eqref{eq:opposed} implies that $\{w_{+1},w_{+2},w_{+3}\}=\{w_{-1},w_{-2},w_{-3}\}$, which contradicts our assumption that the
eigenvalues of $W_+$ and $W_-$ are dif\/ferent.  Consequently $W_+$ and $W_-$ have the same eigenvalues.  This f\/inishes the proof
of Theorem~\ref{thm:slice}:  Any integrable Killing tensor on $S^3$ has a~diagonalisable algebraic curvature tensor.

\subsection{Solution of the second integrability condition}
\label{sub:AICII}

It remains to prove Theorem~\ref{thm:det}.  In the preceding subsections we showed that an algebraic curvature tensor which
satisf\/ies the algebraic integrability conditions is diagonal in a suitable orthogonal basis of~$V$ and that any diagonal
algebraic curvature tensor satisf\/ies the f\/irst of the two algebraic integrability conditions.  We therefore now solve the second
integrability condition~\eqref{eq:AICII} for a diagonal algebraic curvature tensor in the form~\eqref{eq:chiral}.  To this aim
consider the tensor
\begin{gather}
	\label{eq:tensor}
	g^{ij}
	g^{kl}
	R_{ib_1a_2b_2}
	R_{ja_1  kc_1}
	R_{ld_1c_2d_2}
\end{gather}
appearing on the left hand side of \eqref{eq:AICII}.  Suppose it does not vanish.  Then we have $b_1\in\{a_2,b_2\}$ and
$d_1\in\{c_2,d_2\}$, so without loss of generality we can assume $b_1=b_2$ and $d_1=d_2$.  Then only terms with $i=j=a_2$ and
$k=l=c_2$ contribute to the contraction and hence $\{a_1,a_2\}=\{c_1,c_2\}$.  If now $a_2=c_2$, then~\eqref{eq:tensor} vanishes
under complete antisymmetrisation in $a_2$, $b_2$, $c_2$, $d_2$.  This means that the left hand side of \eqref{eq:AICII} vanishes unless
$\{a_1,b_1,c_1,d_1\}=\{a_2,b_2,c_2,d_2\}$.  Note that $a_2$, $b_2$, $c_2$, $d_2$ and therefore $a_1$, $b_1$, $c_1$, $d_1$ have to be pairwise
dif\/ferent due to the antisymmetrisation.

For $\dim V=4$ in particular, \eqref{eq:AICII} reduces to a \emph{sole} condition, which can be written as
\[
	\det
	\begin{pmatrix}
			1&R_{0101}+R_{2323}&R_{0101}R_{2323}\\
			1&R_{0202}+R_{3131}&R_{0202}R_{3131}\\
			1&R_{0303}+R_{1212}&R_{0303}R_{1212}
	\end{pmatrix}
	=0.
\]
Substituting \eqref{eq:chiral} yields
\begin{gather}
	\label{eq:Vandermonde}
	\det
	\begin{pmatrix}
			1&w_1&w_1^2-t_1^2\\
			1&w_2&w_2^2-t_2^2\\
			1&w_3&w_3^2-t_3^2
	\end{pmatrix}
	=0.
\end{gather}
The form of this equation allows us to give an isometry invariant reformulation of the algebraic integrability conditions, which
does not rely on diagonalising the algebraic curvature tensor.  This also entails the isometry invariants in
Proposition~\ref{prop:invariants}.

\begin{Proposition}
	The algebraic integrability conditions for an algebraic curvature tensor with Weyl part~$\mathbf W$ and $($not necessarily
	trace free$)$ Ricci part $\mathbf T$ are equivalent to
	\begin{subequations}
		\begin{gather}
			\label{eq:[W,T]=0}
			[\mathbf W,\mathbf T] =0,\\
			\label{eq:tr(*W^k)=0}
			\tr\big(*\mathbf W^2\big)=
			\tr\big(*\mathbf W^3\big) =0
		\end{gather}
	\end{subequations}
	and the linear dependence of $\mathbf W$, $\mathbf W^2-\mathbf T^2$ and the identity matrix.
\end{Proposition}
\begin{proof}
	This is a consequence of Theorem~\ref{thm:slice} and Proposition~\ref{prop:diagonalisability}, together with
	Corollary~\ref{cor:diagonal->AICI} and the second integrability condition in the form \eqref{eq:Vandermonde}.  Indeed,
	\eqref{eq:[W,T]=0} is the def\/inition of alignedness for $\mathbf R$ and $\tr(*\mathbf W^k)=0$ is equivalent to $\tr
	W_+^k=\tr W_-^k$.  For $k=1$ this equation is the Bianchi identity in the form \eqref{eq:Bianchi}.  Hence
	\eqref{eq:tr(*W^k)=0} means that $W_+$ and $W_-$ have the same characteristic polynomial.  The linear dependence of $\mathbf
	W$, $\mathbf W^2-\mathbf T^2$ and the identity matrix is equivalent to \eqref{eq:Vandermonde}.
\end{proof}

For $t_1=t_2=t_3=0$ the left hand side of \eqref{eq:Vandermonde} is the Vandermonde determinant.  Therefore
\begin{gather}
	\label{eq:quadric}
	(w_1-w_2) t_3^2+(w_2-w_3) t_1^2+(w_3-w_1) t_2^2
	+(w_1-w_2)(w_2-w_3)(w_3-w_1)=0.
\end{gather}
This proves the f\/irst part of Theorem~\ref{thm:det}, for this equation can be written in the form
\[
	\det
	\begin{pmatrix}
		w_2-w_3&-t_3&t_2\\
		t_3&w_3-w_1&-t_1\\
		-t_2&t_1&w_1-w_2
	\end{pmatrix}
	=0
	 .
\]
The second part follows from the residual action of the isometry on diagonal algebraic curvature tensors as described in
Section~\ref{sub:remainder}.

Regarding the eigenvalues $w_1,w_2,w_3\in{\mathbb R}$ as parameters subject to the restriction
\begin{gather}
	\label{eq:tr(W)=0}
	w_1+w_2+w_3=0,
\end{gather}
we can solve this equation for $t_1,t_2,t_3\in{\mathbb R}$, depending on the number of equal eigenvalues:
\begin{enumerate}\itemsep=0pt
	\item
		If $w_1$, $w_2$, $w_3$ are pairwise dif\/ferent with $w_\alpha<w_\beta<w_\gamma$ then \eqref{eq:quadric} describes a one sheeted
		hyperboloid
		\[
			\left(\frac{t_\alpha}{a_\alpha}\right)^2-
			\left(\frac{t_\beta}{a_\beta}\right)^2+
			\left(\frac{t_\gamma}{a_\gamma}\right)^2=1
		\]
		with semi axes
			\begin{gather*}
				a_\alpha =\sqrt{(w_\beta-w_\alpha)(w_\gamma-w_\alpha)},\qquad
				a_\beta =\sqrt{(w_\beta-w_\alpha)(w_\gamma-w_\beta)},\\
				a_\gamma =\sqrt{(w_\gamma-w_\alpha)(w_\gamma-w_\beta)}.
			\end{gather*}
	\item
		If $w_\alpha=w_\beta\not=w_\gamma$, then \eqref{eq:tr(W)=0} implies that the solutions of \eqref{eq:quadric} are those
		$(t_1,t_2,t_3)\in{\mathbb R}^3$ with $t_\beta=\pm t_\alpha$.  Since $(t_\alpha,-t_\beta,-t_\gamma)$ and
		$(t_\alpha,t_\beta,t_\gamma)$ are equivalent under the isometry group, we can assume $t_\beta=t_\alpha$ without loss of
		generality.
	\item
		If $w_1=w_2=w_3$, then \eqref{eq:tr(W)=0} implies that $w_1=w_2=w_3=0$.  In this case \eqref{eq:quadric} is satisf\/ied
		for any choice of $t_1,t_2,t_3\in{\mathbb R}$.  In other words, an algebraic curvature tensor with zero Weyl component is always
		integrable.
\end{enumerate}

\section{The algebraic geometry of the KS-variety}
\label{sec:algebraic}

We have seen that under isometries any integrable Killing tensor is equivalent to one whose algebraic curvature tensor is
diagonal and that the space of such tensors is described by the following variety.
\begin{Definition}
	Let $\mathcal A\isom{\mathbb R}^5$ be the vector space of matrices whose symmetric part is diagonal and trace free, i.e.\ of
	matrices of the form
	\begin{gather}
		\label{eq:matrixx}
		M=
		\begin{pmatrix}
			\Delta_1&-t_3&t_2\\
			t_3&\Delta_2&-t_1\\
			-t_2&t_1&\Delta_3
		\end{pmatrix}
	\end{gather}
	with
	\[
		\tr M=0.
	\]
	We denote by $\AV\subset\mathcal A$ the algebraic variety def\/ined by the equation
	\[
		\det M=0
	\]
	and call its projectivisation $\PV\subset\mathbb P\mathcal A\isom\mathbb P^4$ the {\it Killing--St\"ackel variety
	$($KS-variety$)$}.  A matrix in~$\AV$ will be called a {\it Killing--St\"ackel matrix $($KS-matrix$)$}.
\end{Definition}

\begin{Remark}
	\label{rem:remainder}
	The residual isometry group action on diagonal algebraic curvature tensors def\/ines a natural action of the permutation group
	$S_4$ on the Killing--St\"ackel variety.  This action is given by conjugation with matrices in $\SOG(3)$ under the embedding
	$S_4\subset\SOG(3)$ def\/ined in Section~\ref{sub:remainder}.  Later it will be useful to consider~$S_4\subset\SOG(3)$
	as the symmetry group of an octahedron in~$\mathbb R^3$ with vertices~$\pm e_i$ and adjacent faces opposedly oriented.
\end{Remark}

The Killing--St\"ackel variety $\PV$ is a $3$-dimensional projective subvariety in $\mathbb P^4$ and the quotient~$\PV/S_4$
encodes all information on integrable Killing tensors modulo isometries, the metric and scalar multiples.  In this section we
will investigate its structure from a purely algebraic geometric point of view.  This will f\/inally lead to a complete and
explicit algebraic description of St\"ackel systems.

Recall that St\"ackel systems are $n$-dimensional vector spaces of mutually commuting integrable Killing tensors.  Note that every
St\"ackel system contains the metric and that the metric corresponds to the zero KS-matrix.  Therefore, if we assume that the
algebraic curvature tensors of the Killing tensors in a St\"ackel system are mutually diagonalisable, then St\"ackel systems on~$S^3$ correspond to $2$-dimensional planes in $\AV$ or projective lines in the projectivisation~$\PV$.  We will see a posteriori
that our assumption is justif\/ied.  Hence St\"ackel systems constitute a subvariety of the variety of projective lines on~$\PV$,
also called the Fano variety of $\PV$ and denoted by~$F_1(\PV)$.

Let $\mathcal D_3$ be the generic determinantal variety.  This is the variety of $3\times3$ matrices with vanishing determinant
or, equivalently, matrices of rank one or two.  By def\/inition $\PV$ is a~subvariety of~$\mathbb P\mathcal D_3$.  Hence its Fano
variety~$F_1(\PV)$ is a~subvariety of the Fano variety~$F_1(\mathbb P\mathcal D_3)$.  But the latter is well understood.  It
contains dif\/ferent types of projective linear spaces:
\begin{enumerate}\itemsep=0pt
	\item[1)] projective spaces of matrices whose kernel contains a f\/ixed line,
	\item[2)] projective spaces of matrices whose image is contained in a f\/ixed plane and
	\item[3)] 
the projective plane of antisymmetric matrices.
\end{enumerate}
This motivates why we seek for such spaces in $\PV$.  Note that the projective plane of antisymmetric matrices is obviously
contained in $\PV$.
\begin{Definition}
	We call a projective subspace in the KS-variety consisting of matrices whose kernel contains a f\/ixed line (respectively
	whose image is contained in a f\/ixed plane) an {\it isokernel space} (respectively an {\it isoimage space}).
\end{Definition}
A $3\times3$ matrix $M$ of rank two has a $2$-dimensional image and a $1$-dimensional kernel.  Both are given by its adjugate
matrix $\Adj M$, which is the transpose of the cofactor matrix and satisf\/ies
\begin{gather}
	\label{eq:adjugate}
	(\Adj M)M=M(\Adj M)=(\det M)\Id .
\end{gather}
The adjugate of a matrix $M\in\AV$ is
\begin{gather}
	\label{eq:Adj(M)}
	\Adj M=
	\left(
		\begin{matrix}
			t_1t_1 + \Delta_2\Delta_3& t_2t_1 + \Delta_3     t_3 &t_3t_1 - \Delta_2     t_2\\
			t_1t_2  - \Delta_3     t_3& t_2t_2 + \Delta_3\Delta_1 &t_3t_2 + \Delta_1     t_1\\
			t_1t_3 + \Delta_2     t_2& t_2t_3 - \Delta_1     t_1 &t_3t_3 + \Delta_1\Delta_2
		\end{matrix}
	\right) .
\end{gather}
The adjugate matrix of a rank two $3\times3$ matrix has rank one, which means that the columns of~$\Adj M$ span a line in~${\mathbb R}^3$.  And because $\det M=0$, we deduce from \eqref{eq:adjugate} that the kernel of~$M$ is the column space of~$\Adj M$.

The adjugate matrix of a rank one $3\times3$ matrix is zero.  It is known, that the singular locus of~$\mathbb P\mathcal D_3$ is
the subvariety of rank one matrices.  The following proposition characterises the singular locus of the KS-variety.  For
simplicity we will be a bit imprecise and not strictly distinguish between a matrix $M\in\AV\setminus\{0\}$ and its image ${\mathbb R}
M\in\PV$ under the canonical projection $\AV\setminus\{0\}\to\PV$.  Depending on the context, we will refer to an element of
$\PV$ as ``matrix'' or ``point''.
\begin{Proposition}[singularities of the KS-variety]
	\label{prop:singular}
	The KS-variety $\PV$ has ten singular points, given by the six rank one matrices
	\begin{subequations}
		\label{eq:singular}
		\begin{gather}
			\label{eq:vertices}
			V_{\pm1} \coloneq
			\begin{pmatrix}
				0&0&0\\
				0&+1&\mp1\\
				0&\pm1&-1
			\end{pmatrix}, \qquad
			V_{\pm2} \coloneq
			\begin{pmatrix}
				-1&0&\pm1\\
				0&0&0\\
				\mp1&0&+1
			\end{pmatrix},\qquad
			V_{\pm3} \coloneq
			\begin{pmatrix}
				+1&\mp1&0\\
				\pm1&-1&0\\
				0&0&0
			\end{pmatrix}
		\end{gather}
		and the four skew symmetric matrices
		\begin{gather}
			\label{eq:centres}
			C_0      \coloneq V_{+1     }+V_{+2    }+V_{+3     },\qquad
			C_\alpha \coloneq V_{-\alpha}+V_{+\beta}+V_{+\gamma},
		\end{gather}
	\end{subequations}
	where $(\alpha,\beta,\gamma)$ denotes an even permutation of $(1,2,3)$.  Moreover, the rank one matrices in $\PV$ are
	exactly those six in \eqref{eq:vertices}.
\end{Proposition}
\begin{Definition}
	We will call the six singularities $V_{\pm\alpha}$ the {\it rank one singular points} in $\PV$ and the four singularities
	$C_i$ the {\it skew symmetric singular points} in $\PV$.
\end{Definition}
\begin{proof}
	$\PV$ is singular at ${\mathbb R} M$ if and only if $\AV$ is singular at $M$.  Since $\AV$ is the zero locus of the determinant
	function $\det\colon\mathcal A\to{\mathbb R}$, this is the case if and only if the derivative of the determinant function at
	$M\in\AV$, given by
	\begin{gather}
		\label{eq:ddet}
		\begin{array}{@{}rrcl}
			(d\det)_M\colon&\mathcal A&\longrightarrow&{\mathbb R},\\
			&A&\mapsto&\tr(A\Adj M),
		\end{array}
	\end{gather}
	is the zero map.  The condition that $\tr(A\Adj M)=0$ for all $A\in\mathcal A$ is equivalent to
	\begin{gather*}
		t_1^2+\Delta_2\Delta_3=
		t_2^2+\Delta_3\Delta_1 =
		t_3^2+\Delta_1\Delta_2,\qquad
		\Delta_1t_1=
		\Delta_2t_2=
		\Delta_3t_3 =0 .
	\end{gather*}
	We leave it to the reader to verify that the solutions $M\in\AV$ of these equations are exactly the matrices
	\eqref{eq:singular} and their multiples.

	The last statement now follows from \eqref{eq:ddet} and the equivalence of $\operatorname{rank}M=1$ and $\Adj M=0$ for
	$3\times3$-matrices.
\end{proof}

\begin{Remark}
	Under the natural $S_4$-action on $\PV$, all rank one singular points are equivalent.  This means that the singularities
	$V_{\pm\alpha}$ are all mapped to a single point~$V$ in the quotient~$\PV/S_4$.  The same holds for the skew symmetric
	singular points, which are mapped to a single point~$C\in\PV/S_4$.
\end{Remark}

We now compute the isokernel spaces in $\PV$.   To f\/ind all matrices $M\in\AV$ annihilating a given vector $\vec n\in{\mathbb R}^3$,
we consider the equation $M\vec n=0$ as an equation in $M$ for a f\/ixed $\vec n=(n_1,n_2,n_3)$, where $M$ is of the form
$M=\Delta+T$ with
\begin{gather*}
	\Delta =
	\begin{pmatrix}
		\Delta_1&0&0\\
		0&\Delta_2&0\\
		0&0&\Delta_3
	\end{pmatrix},\qquad
	T =
	\begin{pmatrix}
		0&-t_3&t_2\\
		t_3&0&-t_1\\
		-t_2&t_1&0
	\end{pmatrix}	 .
\end{gather*}
We can regard $\vec t=(t_1,t_2,t_3)\in{\mathbb R}^3$ as a parameter and solve $\Delta\vec n=-T\vec n$ for $\Delta$.  For the time being,
let us assume $n_1n_2n_3\not=0$.  We then obtain a linear family of matrices
\begin{subequations}
	\label{eq:isokernel}
	\begin{gather}
		\label{eq:isokernel:parametrisation}
		M=
		\begin{pmatrix}
			\frac{n_2}{n_1}t_3-\frac{n_3}{n_1}t_2&-t_3&t_2\\
			t_3&\frac{n_3}{n_2}t_1-\frac{n_1}{n_2}t_3&-t_1\\
			-t_2&t_1&\frac{n_1}{n_3}t_2-\frac{n_2}{n_3}t_1\\
		\end{pmatrix}
		 .
	\end{gather}
	All these matrices satisfy $\det M=0$ and the condition $\tr M=0$ imposes the restriction
	\begin{gather}
		\label{eq:isokernel:restriction}
		\vec t\perp
		\left(
			\frac{n_3}{n_2}-\frac{n_2}{n_3},
			\frac{n_1}{n_3}-\frac{n_3}{n_1},
			\frac{n_2}{n_1}-\frac{n_1}{n_2}
		\right)
	\end{gather}
\end{subequations}
on the parameter $\vec t\in{\mathbb R}^3$.  This shows that the matrices in $\AV$ annihilating a f\/ixed vector $\vec n\not=0$ form a
linear subspace.  If the right hand side of \eqref{eq:isokernel:restriction} is zero, i.e.\ if the condition is void, this
subspace has dimension three and def\/ines a projective plane in $\PV$.  This happens if and only if
\begin{gather}
	\label{eq:exceptional}
	\lvert n_1\rvert=\lvert n_2\rvert=\lvert n_3\rvert.
\end{gather}
If the right hand side of~\eqref{eq:isokernel:restriction} is not zero, this subspace has dimension two and def\/ines a projective
line in $\PV$.  One can check that this also holds true for the case $n_1n_2n_3=0$.  We have shown the following:
\begin{Lemma}
	\label{lem:isokernel}
	If $\vec n\in{\mathbb R}^3$ satisfies \eqref{eq:exceptional}, then the set of matrices in $\PV$ annihilating $\vec n$ is a~projective
	plane of the form~\eqref{eq:isokernel:parametrisation}.  Otherwise it is a projective line.
\end{Lemma}

Recall that the kernel of a rank two matrix of the form \eqref{eq:matrixx} is the row space of the adjugate matrix~\eqref{eq:Adj(M)}.  This def\/ines a rational map
\[
	\pi\colon \ \PV\dashrightarrow\mathbb P^2,
\]
whose f\/ibres are the maximal isokernel spaces.  This map is well def\/ined except for the six rank one matrices in~$\PV$.

We want to give a parametrisation of the isokernel lines that is uniform and more geometric than \eqref{eq:isokernel}.  For this
purpose we def\/ine two embeddings $\iota,\nu\colon\mathbb P^2\to\PV$ given by
\begin{gather}
	\label{eq:embeddings}
	\iota(n) =
	\begin{pmatrix}
		0&-n_3&n_2\\
		n_3&0&-n_1\\
		-n_2&n_1&0
	\end{pmatrix},\qquad
	\nu(n) =
	\begin{pmatrix}
		n_2^2-n_3^2&-n_1n_2&n_3n_1\\
		n_1n_2&n_3^2-n_1^2&-n_2n_3\\
		-n_3n_1&n_2n_3&n_1^2-n_2^2
	\end{pmatrix}
\end{gather}
for $n=(n_1:n_2:n_3)\in\mathbb P^2$.  One readily checks that $\iota(n)$ and $\nu(n)$ indeed lie in $\PV$.  Note that the image
of $\iota$ is the projective plane of skew symmetric matrices in~$\PV$.  Under the inclusion $\PV\subset\mathbb P\mathcal
A\isom\mathbb P^4$ the map $\iota$ is a linear embedding of $\mathbb P^2$ in $\mathbb P^4$ whilst $\nu$ is the composition of
the Veronese embedding
\[
	\begin{array}{@{}rcl}
		\mathbb P^2&\hookrightarrow&\mathbb P^5,\\
		(n_1:n_2:n_3)&\mapsto&\bigl( n_1^2\,:\,n_2^2\,:\,n_3^2\,:\,n_2n_3\,:\,n_3n_1\,:n_1n_2 \bigr),
	\end{array}
\]
and a projection $\mathbb P^5\twoheadrightarrow\mathbb P^4$ under which the Veronese surface remains smooth.

\begin{Proposition}[isokernel spaces in the KS-variety]
	\label{prop:isokernel}
	Recall from Proposition~{\rm \ref{prop:singular}} that the singular points in the KS-variety $\PV$ are the four skew
	symmetric singular points $C_0$, $C_1$, $C_2$, $C_3$ as well as six rank one singular points $V_{\pm\alpha}$, $\alpha=1,2,3$.
	\begin{enumerate}\itemsep=0pt
		\item[$1.$] The projective plane through $V_{+1}$, $V_{+2}$ and $V_{+3}$ is an isokernel plane in $\PV$ and contains~$C_0$.  In
			the same way, the three points $V_{-\alpha}$, $V_{+\beta}$ and $V_{+\gamma}$ define an isokernel plane through
			$C_\alpha$ for $\alpha=1,2,3$.
	\end{enumerate}
	Let the maps $\pi\colon\PV\dashrightarrow\mathbb P^2$ and $\iota,\nu\colon\mathbb P^2\to\PV$ be defined as above.
	\begin{enumerate} 
		\item[$2.$] For any non-singular point $M\in\PV$, the points $M$, $\iota(\pi(M))$ and $\nu(\pi(M))$ are collinear but do not all
			three coincide.  The projective line they define lies in $\PV$ and is an isokernel space through~$M$.  Each of these
			lines contains a unique skew symmetric point, namely $\iota(\pi(M))$. 
	\end{enumerate}
	Moreover, any isokernel space in $\PV$ of maximal dimension is of either of the above forms and hence contains a unique skew
	symmetric point.
\end{Proposition}

\begin{Remark}
	A corresponding characterisation of isoimage spaces follows from the fact that matrix transposition def\/ines an involution of
	$\PV$ which interchanges isokernel and isoimage spaces.
\end{Remark}

\begin{proof}
1~follows directly from Lemma~\ref{lem:isokernel} and the explicit form \eqref{eq:singular} of
	the singular points.  For~2 
notice that $n\coloneq\pi(M)$ is the kernel of $M$ and well def\/ined for
	$M\not=V_\alpha$.  The def\/initions~\eqref{eq:embeddings} of $\iota$ and $\nu$ then show that $n$ is also the kernel of
	$\iota(n)$ and $\nu(n)$.  These def\/initions also show that~$M$,~$\iota(n)$ and~$\nu(n)$ do not all coincide for $M\not=C_i$.
	Their collinearity follows from Lemma~\ref{lem:isokernel} in the case where $n$ does not satisfy~\eqref{eq:exceptional}
	and from $\iota(n)=\nu(n)$ in case it does.  The  last statement is now a consequence of Lemma~\ref{lem:isokernel}.
\end{proof}

\begin{figure}[h]
	\centering
	\includegraphics[width=.55\textwidth]{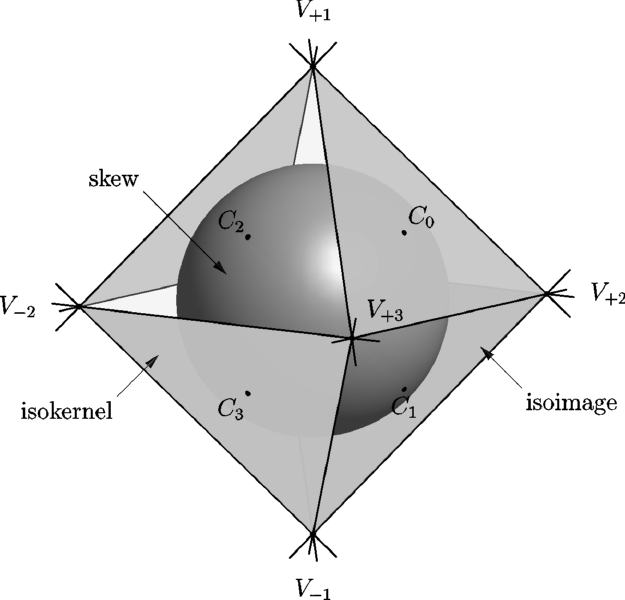}
	\caption{Singularities and projective planes in the KS-variety.}
	\label{fig:octahedron}
\end{figure}

We can illustrate the content of Proposition~\ref{prop:isokernel} geometrically as follows.  This is depicted in
Fig.~\ref{fig:octahedron}.  The six rank one singular points of $\PV\subset\mathbb P^4$ constitute the vertices of a
regular octahedron in $\mathbb P^4$ whose faces (and the planes they def\/ine) are entirely contained in $\PV$.  The set of eight
faces is divided into two sets of non-adjacent faces, corresponding to the four isokernel planes (shaded) respectively the four
isoimage planes in $\PV$.  Opposite faces intersect in their respective centres.  These are the skew symmetric singular points
of $\PV$ and the intersection points of the octahedron with the projective plane of skew symmetric matrices.  In
Fig.~\ref{fig:octahedron} this ninth projective plane in $\PV$ is depicted as the insphere of the octahedron.

The permutation group $S_4$ acts on this conf\/iguration by symmetries of the octahedron preserving the two sets of non-adjacent
faces.  Matrix transposition corresponds to a point ref\/lection, which exchanges opposite faces and completes the $S_4$-action to
the full octahedral symmetry.

In anticipation of the results in the next section we introduce the following name:
\begin{Definition}
	\label{def:Staeckel-lines}
	We call the isokernel lines of the form~2
in Proposition~\ref{prop:isokernel} {\it ``St\"ackel
	lines''}.
\end{Definition}
Every St\"ackel line intersects the projective plane of skew symmetric matrices transversely.  Conversely, every point on this
plane determines a unique St\"ackel line unless it is one of the four singular skew symmetric points.  The St\"ackel lines through
these points are exactly those projective lines in each isokernel face which pass through the respective face center.

\section{St\"ackel systems}
\label{sec:Staeckel}

In this section we express the condition that two Killing tensors commute as an algebraic condition on the corresponding
algebraic curvature tensors.  This will eventually justify the term ``St\"ackel lines'' chosen in
Def\/inition~\ref{def:Staeckel-lines} by proving that these lines indeed correspond to St\"ackel systems.

\begin{Remark}
	Most results in this and the following section are valid for an arbitrary non-f\/lat constant curvature (pseudo-)Riemannian
	manifold $M$ that is embedded into some vector space~$V$ as the hypersurface
	\[
		M=\{v\in V\colon g(v,v)=1\}\subset V
	\]
	of length one vectors with respect to a non-degenerated symmetric bilinear form $g$ on $V$.  Note that any non-f\/lat constant
	curvature (psudo-)Riemannian manifold is locally isometric to such a model.
\end{Remark}

\subsection{Young tableaux}

Throughout this section we use Young tableaux in order to characterise index symmetries of tensors.  Young tableaux def\/ine
elements in the group algebra of the permutation group~$S_d$.  That is, a Young tableau stands for a (formal) linear combination
of permutations of $d$ objects.  In our case, these objects will be certain tensor indices.  For the sake of simplicity of
notation we will identify a Young tableau with the group algebra element it def\/ines.  It will be suf\/f\/icient for our purposes to
consider the examples below.  For more examples and a detailed explanation of the notation as well as the techniques used here,
we refer the reader to \cite{Schoebel}.

A Young tableau consisting of a single row denotes the sum of all permutations of the indices in this row.  For example, using
cycle notation,
\[
	\young(\atwo\cone\ctwo)=e+(a_2c_1)+(c_1c_2)+(c_2a_2)+(a_2c_1c_2)+(c_2c_1a_2).
\]
This is an element in the group algebra of the group of permutations of the indices~$a_2$,~$c_1$ and~$c_2$ (or any superset).
In the same way a Young tableau consisting of a single column denotes the signed sum of all permutations of the indices in this
column, the sign being the sign of the permutation.  For example,
\[
	\young(\aone,\bone,\dtwo)=e-(a_1b_1)-(b_1d_2)-(d_2a_1)+(a_1b_1d_2)+(d_2b_1a_1) .
\]
We call these {\it row symmetrisers} respectively {\it column antisymmetrisers}.  The reason we def\/ine them without the usual
normalisation factors is that then all numerical constants appear expli\-cit\-ly in our computations (although irrelevant for our
concerns).

The group multiplication extends linearly to a natural product in the group algebra.  A~general Young tableau is then simply the
product of all column antisymmetrisers and all row symmet\-ri\-sers of the tableau.  We will only deal with Young tableaux having a
``hook shape'', such as the following:
\begin{gather}
	\label{eq:hook}
	\young(\aone\atwo\cone\ctwo,\bone,\dtwo)=\young(\aone\atwo\cone\ctwo)\young(\aone,\bone,\dtwo) .
\end{gather}

The inversion of group elements extends linearly to an involution of the group algebra.  If we consider elements in the group
algebra as linear operators on the group algebra itself, this involution is the adjoint with respect to the natural inner
product on the group algebra, given by def\/ining the group elements to be an orthonormal basis.  Since this operation leaves
symmetrisers and antisymmetrisers invariant, it simply exchanges the order of symmetrisers and antisymmetrisers in a Young
tableau.  The adjoint of~\eqref{eq:hook} for example is
\begin{gather}
	\label{eq:antihook}
	{\young(\aone\atwo\cone\ctwo,\bone,\dtwo)}^\adjoint=\young(\aone,\bone,\dtwo)\young(\aone\atwo\cone\ctwo).
\end{gather}
Properly scaled, Young tableaux def\/ine projectors onto irreducible $S_d$-representations.  We are now able to state a lemma
which expresses the decomposition
\[
	\young(~,~)\otimes\young(~~~~)
	\isom
	\young(~~~~,~,~)
	\oplus
	\young(~~~~~,~)
	 ,
\]
given by the Littlewood--Richardson rule, in terms of Young projectors.  Its proof can be found in \cite{Schoebel}.
\begin{Lemma}
	\label{lem:hooks}
	In the group algebra of the permutation group of the indices $a_1$, $a_2$, $c_1$, $c_2$, $b_1$ and~$d_2$ the following
	identity holds:
	\begin{gather}
		\label{eq:projectors}
		\frac1{2!}\young(\bone,\dtwo)\cdot\frac1{4!}\young(\aone\atwo\cone\ctwo)
		=
		\frac1{10368}{\young(\aone\atwo\cone\ctwo,\bone,\dtwo)}         {\young(\aone\atwo\cone\ctwo,\bone,\dtwo)}^\adjoint+
		\frac1{34560}{\young(\bone\aone\atwo\cone \ctwo,\dtwo)}^\adjoint{\young(\bone\aone\atwo\cone \ctwo,\dtwo)}
		 .
	\end{gather}
\end{Lemma}

The action of the permutation group on tensors extends linearly to an action of the group algebra.  In particular, any Young
tableau acts on tensors with corresponding indices.  For example,
\[
	\young(\bone,\atwo,\ctwo)
	 T_{b_1a_2c_2}
	=T_{b_1a_2c_2}
	-T_{a_2b_1c_2}
	-T_{b_1c_2a_2}
	-T_{c_1a_2b_2}
	+T_{a_1c_2b_2}
	+T_{c_1b_2a_2}
	 .
\]
To give another example, the operator \eqref{eq:hook} acts on a tensor $T_{b_1b_2d_1d_2a_2c_2}$ by an antisymmetrisation in the
indices $b_1$, $a_2$, $c_2$ and a subsequent symmetrisation in the indices $b_1$, $b_2$, $d_1$, $d_2$.  In the same way its
adjoint \eqref{eq:antihook} acts by f\/irst symmetrising and then antisymmetrising.

\subsection{Commuting Killing tensors}

In the following we express the fact that two Killing tensors commute (as endomorphisms) as a~purely algebraic condition on
their corresponding algebraic curvature tensors.  Notice that the commuting of Killing tensors as endomorphisms on the tangent
space is not equivalent to the commuting of their respective algebraic curvature tensors as endomorphisms on~$\Lambda^2V$, as
one may guess na\"{\i}vely.
\begin{Proposition}[algebraic counterpart of the commutator] \label{prop:commutator}
	Let $K$ and $\tilde K$ be two Killing tensors on a non-flat constant curvature manifold with algebraic
	curvature tensors $R$ and $\tilde R$.  Then the following statements are equivalent
	\begin{subequations}
		\begin{gather}
			\label{eq:commutator}
			[K,\tilde K] =0,\\
			\label{eq:commutator:hook}
			{\young(\bone\btwo\done\dtwo,\atwo,\ctwo)}^\adjoint
			g_{ij}R\indices{^i_{b_1a_2b_2}}\tilde R\indices{^j_{d_1c_2d_2}}
			 =0,\\[\medskipamount]
			\label{eq:commutator:4+2-}
			{\young(\atwo,\ctwo)}{\young(\bone\btwo\done\dtwo)}
			g_{ij}R\indices{^i_{b_1a_2b_2}}\tilde R\indices{^j_{d_1c_2d_2}}
			 =0.
		\end{gather}
	\end{subequations}
\end{Proposition}
\begin{proof}
	Using $\nabla_vx^b=v^b$ we write \eqref{eq:correspondence} in local coordinates for coordinate vectors $v=\partial_\alpha$
	and $w=\partial_\beta$:
	\[
		K_{\alpha\beta}
		=R_{a_1b_1a_2b_2}x^{a_1}x^{a_2}\nabla_\alpha x^{b_1}\nabla_\beta x^{b_2}
		 .
	\]
	The product of $K$ and $\tilde K$, regarded as endomorphisms, is then given by
	\[
		K\indices{^\alpha_\gamma}\tilde K\indices{^\gamma_\beta}
		=
		R_{a_1b_1a_2b_2}\tilde R_{c_1d_1c_2d_2}
		x^{a_1}x^{a_2}
		x^{c_1}x^{c_2}
		\nabla^\alpha x^{b_1}\nabla_\gamma x^{b_2}
		\nabla^\gamma x^{d_1}\nabla_\beta  x^{d_2} .
	\]
	In \cite{Schoebel} we proved the following identity:
	\[
		\nabla_\gamma x^{b_2}\nabla^\gamma x^{d_1}=g^{b_2d_1}-x^{b_2}x^{d_1} .
	\]
	As a consequence of the antisymmetry of algebraic curvature tensors in the last index pair, the term $x^{b_2}x^{d_1}$ does
	not contribute when substituting this identity into the previous expression:
	\[
		K\indices{^\alpha_\gamma}\tilde K\indices{^\gamma_\beta}
		=g^{b_2d_1}R_{a_1b_1a_2b_2}\tilde R_{c_1d_1c_2d_2}x^{a_1}x^{a_2}x^{c_1}x^{c_2}\nabla^\alpha x^{b_1}\nabla_\beta x^{d_2} .
	\]
	The commutator $[K,\tilde K]$ is therefore given by
	\[
		[K,\tilde K]_{\alpha\beta}
		=g^{b_2d_1}R_{a_1b_1a_2b_2}\tilde R_{c_1d_1c_2d_2}x^{a_1}x^{a_2}x^{c_1}x^{c_2}\nabla_{[\alpha}x^{b_1}\nabla_{\beta]}x^{d_2}
	\]
	and vanishes if and only if
	\begin{gather}
		\label{eq:intermediate}
		g^{b_2d_1}R_{a_1b_1a_2b_2}\tilde R_{c_1d_1c_2d_2}x^{a_1}x^{a_2}x^{c_1}x^{c_2}v^{[b_1}w^{d_2]}=0
	\end{gather}
	for all $x\in M$ and $v,w\in T_xM$.  That is, for all $x,v,w\in V$ with
	\begin{gather}
		\label{eq:restrictions}
		g(x,x) =1, \qquad
		g(x,v)=g(x,w) =0.
	\end{gather}
	We can drop the restriction $g(x,x)=1$ since ${\mathbb R} M\subseteq V$ is open.  We can also drop the restrictions $g(x,v)=0$ and
	$g(x,w)=0$ by decomposing arbitrary vectors $v,w\in V$ under the decomposition $V=T_xM\oplus{\mathbb R} x$.  To see this, notice that
	\eqref{eq:intermediate} is trivially satisf\/ied for $v=x$ or for $w=x$.  Indeed, in this case the tensor
	\begin{gather}
		\label{eq:RR}
		R_{a_1b_1a_2b_2}\tilde R_{c_1d_1c_2d_2}
	\end{gather}
	is implicitly symmetrised over f\/ive indices and Dirichlet's drawer principle tells us that this comprises a symmetrisation
	in one of the four antisymmetric index pairs.  This means we can omit the restrictions \eqref{eq:restrictions} completely.
	In other words, $[K,\tilde K]=0$ is equivalent to \eqref{eq:intermediate} being satisf\/ied for \emph{all} $x,v,w\in V$.

	Now notice that the tensor
	\[
		x^{a_1}x^{a_2}x^{c_1}x^{c_2}v^{[b_1}w^{d_2]}
	\]
	appearing in \eqref{eq:intermediate} is completely symmetric in the indices $a_1$, $a_2$, $c_1$ and $c_2$ and completely
	antisymmetric in the indices~$b_1$ and~$d_2$.  Applying \eqref{eq:projectors} to this tensor therefore yields
	\begin{gather}
			x^{a_1}x^{a_2}x^{c_1}x^{c_2}v^{[b_1}w^{d_2]}
			 =\frac1{2!}\young(\bone,\dtwo)\cdot\frac1{4!}\young(\aone\atwo\cone\ctwo)x^{a_1}x^{a_2}x^{c_1}x^{c_2}v^{[b_1}w^{d_2]}\nonumber\\[\medskipamount]
\hphantom{x^{a_1}x^{a_2}x^{c_1}x^{c_2}v^{[b_1}w^{d_2]}}{}		
=\frac1{10368}\cdot{\young(\aone\atwo\cone\ctwo,\bone,\dtwo)}{\young(\aone\atwo\cone\ctwo,\bone,\dtwo)}^\adjoint x^{a_1}x^{a_2}x^{c_1}x^{c_2}v^{[b_1}w^{d_2]}\nonumber\\[\medskipamount]
\hphantom{x^{a_1}x^{a_2}x^{c_1}x^{c_2}v^{[b_1}w^{d_2]}}{}	
+\frac1{34560}\cdot{\young(\bone\aone\atwo\cone \ctwo,\dtwo)}^\adjoint{\young(\bone\aone\atwo\cone\ctwo,\dtwo)}x^{a_1}x^{a_2}x^{c_1}x^{c_2}v^{[b_1}w^{d_2]}.
\label{eq:decomposition}
	\end{gather}
	From the def\/inition of the Young tableaux and their adjoints,
	\begin{alignat*}{3}
	&	{\young(\aone\atwo\cone\ctwo,\bone,\dtwo)}          ={\young(\aone\atwo\cone\ctwo)}{\young(\aone,\bone,\dtwo)},\qquad &&
		{\young(\aone\atwo\cone\ctwo,\bone,\dtwo)}^\adjoint ={\young(\aone,\bone,\dtwo)}{\young(\aone\atwo\cone\ctwo)}, &\\
	&	{\young(\bone\aone\atwo\cone \ctwo,\dtwo)}^\adjoint ={\young(\bone,\dtwo)}{\young(\bone\aone\atwo\cone\ctwo)}, \qquad &&
		{\young(\bone\aone\atwo\cone\ctwo,\dtwo)}           ={\young(\bone\aone\atwo\cone\ctwo)}{\young(\bone,\dtwo)},&
	\end{alignat*}
	together with the properties
	\begin{gather*}
		{\young(\aone\atwo\cone\ctwo)}x^{a_1}x^{a_2}x^{c_1}x^{c_2} =4!x^{a_1}x^{a_2}x^{c_1}x^{c_2},\qquad
		{\young(\bone,\dtwo)}v^{[b_1}w^{d_2]} =2{\young(\bone,\dtwo)}v^{b_1}w^{d_2}
	\end{gather*}
	and
	\begin{gather*}
		{\young(\aone,\bone,\dtwo)}^2 =3!{\young(\aone,\bone,\dtwo)},\qquad
		{\young(\bone\aone\atwo\cone\ctwo)}^2 =5!{\young(\bone\aone\atwo\cone\ctwo)}
	\end{gather*}
	we see that \eqref{eq:decomposition} simplif\/ies to
	\begin{gather*}
		x^{a_1}x^{a_2}x^{c_1}x^{c_2}v^{[b_1}w^{d_2]}
		=
		\left(
			\frac{4!\cdot3!}{10368}{\young(\aone\atwo\cone\ctwo,\bone,\dtwo)}+
			\frac{2 \cdot5!}{34560}{\young(\bone\aone\atwo\cone \ctwo,\dtwo)}^\adjoint
		\right)
		x^{a_1}x^{a_2}x^{c_1}x^{c_2}v^{b_1}w^{d_2}.
	\end{gather*}
	When substituted into \eqref{eq:intermediate}, this yields
	\begin{gather}
		\label{eq:upper}
		g^{b_2d_1}R_{a_1b_1a_2b_2}\tilde R_{c_1d_1c_2d_2}
		\left(
			\frac1{72} {\young(\aone\atwo\cone\ctwo,\bone,\dtwo)}+
			\frac1{144}{\young(\bone\aone\atwo\cone \ctwo,\dtwo)}^\adjoint
		\right)
		x^{a_1}x^{a_2}x^{c_1}x^{c_2}v^{b_1}w^{d_2}
		=0		 .
	\end{gather}
	As elements of the group algebra of the permutation group, the Young tableaux are linear combinations of permutations
	$\pi$ of the indices $a_1$, $a_2$, $c_1$, $c_2$, $b_1$, $d_2$.  For any such permutation $\pi$ we have
	\begin{gather*}
		g^{b_2d_1}R_{a_1b_1a_2b_2}\tilde R_{c_1d_1c_2d_2}
		x^{\pi(a_1)}x^{\pi(a_2)}
		x^{\pi(c_1)}x^{\pi(c_2)}
		v^{\pi(b_1)}w^{\pi(d_2)}\\
		\qquad{}=
		g^{b_2d_1}R_{\pi^\adjoint(a_1)\pi^\adjoint(b_1)\pi^\adjoint(a_2)b_2}\tilde R_{\pi^\adjoint(c_1)d_1\pi^\adjoint(c_2)\pi^\adjoint(d_2)}
		x^{a_1}x^{a_2}
		x^{c_1}x^{c_2}
		v^{b_1}w^{d_2},
	\end{gather*}
	where $\pi^\adjoint=\pi^{-1}$ is the adjoint of $\pi$.  We can thus replace the Young tableaux in \eqref{eq:upper}
	acting on upper indices by its adjoint acting on lower indices:
	\begin{gather}
		\left(
			\left(\!
				\frac1{72} {\young(\aone\atwo\cone\ctwo,\bone,\dtwo)}^\adjoint+
				\frac1{144}{\young(\bone\aone\atwo\cone \ctwo,\dtwo)}
			\!\right)\!
		g^{b_2d_1}R_{a_1b_1a_2b_2}\tilde R_{c_1d_1c_2d_2}
		\right)
		x^{a_1}x^{a_2}x^{c_1}x^{c_2}v^{b_1}w^{d_2}
		=0.\!\!\!\label{eq:lower}
	\end{gather}
	Now notice that the second Young tableau involves a symmetrisation over the f\/ive indices $b_1$, $a_1$, $a_2$, $c_1$, $c_2$ and that, as
	above, the symmetrisation of the tensor \eqref{eq:RR} in any f\/ive indices is zero.  Hence the second term in~\eqref{eq:lower} vanishes and we obtain
	\[
		\left(
			{\young(\aone\atwo\cone\ctwo,\bone,\dtwo)}^\adjoint
			g^{b_2d_1}R_{a_1b_1a_2b_2}\tilde R_{c_1d_1c_2d_2}
		\right)
		x^{a_1}x^{a_2}x^{c_1}x^{c_2}v^{b_1}w^{d_2}
		=0 .
	\]
	Recall that $[K,\tilde K]=0$ is equivalent to this condition being satisf\/ied for all $x,v,w\in V$.  By polarising in $x$ we
	get
	\[
		\left(
			{\young(\aone\atwo\cone\ctwo,\bone,\dtwo)}^\adjoint
			g^{b_2d_1}R_{a_1b_1a_2b_2}\tilde R_{c_1d_1c_2d_2}
		\right)
		x^{a_1}y^{a_2}
		z^{c_1}t^{c_2}
		v^{b_1}w^{d_2}
		=0
	\]
	for all $x,y,z,t,v,w\in V$ and hence
	\[
		{\young(\aone\atwo\cone\ctwo,\bone,\dtwo)}^\adjoint
		g^{b_2d_1}R_{a_1b_1a_2b_2}\tilde R_{c_1d_1c_2d_2}
		=0.
	\]
	This is the same as \eqref{eq:commutator:hook} after appropriately renaming, lowering and rising indices.  We have proven
	the equivalence \eqref{eq:commutator} $\Leftrightarrow$ \eqref{eq:commutator:hook}.

	We now prove the equivalence \eqref{eq:commutator:hook} $\Leftrightarrow$ \eqref{eq:commutator:4+2-}.  Start from
	\eqref{eq:commutator:hook} by expanding the Young tableau:
	\[
		{\young(\bone,\atwo,\ctwo)}
		{\young(\bone\btwo\done\dtwo)}
		g_{ij}R\indices{^i_{b_1a_2b_2}}\tilde R\indices{^j_{d_1c_2d_2}}
		=0
		 .
	\]
	In order to sum over all $4!$ permutations when carrying out the symmetrisation in the indices $b_1$, $b_2$, $d_1$, $d_2$, one can
	f\/irst take the sum over the $4$ cyclic permutations of $b_1$, $b_2$, $d_1$, $d_2$, then f\/ix the index $b_1$ and f\/inally sum over all
	$3!$ permutations of the remaining $3$ indices $b_2$, $d_1$, $d_2$:
	\begin{gather*}
			{\young(\bone,\atwo,\ctwo)}
			{\young(\btwo\done\dtwo)}
			g_{ij}\!
			\Bigl(
				 R\indices{^i_{\underline b_1\underline a_2           b_2}}\tilde R\indices{^j_{           d_1\underline c_2           d_2}}\!+
				 R\indices{^i_{           b_2\underline a_2           d_1}}\tilde R\indices{^j_{           d_2\underline c_2\underline b_1}}\! +
				 R\indices{^i_{           d_1\underline a_2           d_2}}\tilde R\indices{^j_{\underline b_1\underline c_2           b_2}}\!+
				 R\indices{^i_{           d_2\underline a_2\underline b_1}}\tilde R\indices{^j_{           b_2\underline c_2           d_1}}
			\Bigr)
			=0.
	\end{gather*}
	For a better readability we underlined each antisymmetrised index.  Permuting the indices of the terms in the parenthesis
	under symmetrisation in $b_2$, $d_1$, $d_2$ and antisymmetrisation in $b_1$, $a_2$, $c_2$ we can gather the f\/irst and last as well as
	the second and third term:
	\begin{gather*}
			{\young(\bone,\atwo,\ctwo)}
			{\young(\btwo\done\dtwo)}
			g_{ij}
			\Bigl(
				 (R\indices{^i_{\underline b_1\underline a_2b_2}}+R\indices{^i_{b_2\underline a_2\underline b_1}})\tilde R\indices{^j_{d_1\underline c_2d_2}} +
				 R\indices{^i_{d_1\underline a_2d_2}}(\tilde R\indices{^j_{\underline b_1\underline c_2b_2}}+\tilde R\indices{^j_{b_2\underline c_2\underline b_1}})
			\Bigr)
			=0.
	\end{gather*}
	Using the symmetries of $R_{b_1a_2b_2}$, the terms in the inner parentheses can be rewritten as
	\begin{gather*}
		R\indices{^i_{\underline b_1\underline a_2b_2}}
		 =-R\indices{^i_{\underline b_1b_2\underline a_2}},\qquad
		R\indices{^i_{b_2\underline a_2\underline b_1}}
		=-R\indices{^i_{\underline a_2\underline b_1b_2}}-R\indices{^i_{\underline b_1b_2\underline a_2}}
		= R\indices{^i_{\underline a_2b_2\underline b_1}}-R\indices{^i_{\underline b_1b_2\underline a_2}}
	\end{gather*}
	resulting in
	\[
		{\young(\bone,\atwo,\ctwo)}
		{\young(\btwo\done\dtwo)}
		g_{ij}
		\Bigl(
			R\indices{^i_{\underline a_2           b_2\underline b_1}}\tilde R\indices{^j_{           d_1\underline c_2           d_2}}+
			R\indices{^i_{           d_1\underline a_2           d_2}}\tilde R\indices{^j_{\underline c_2           b_2\underline b_1}}
		\Bigr)
		=0		 .
	\]
	As above, when carrying out the antisymmetrisation over $b_1$, $a_2$, $c_2$, we can f\/irst sum over the three cyclic permutations
	of $b_1$, $a_2$, $c_2$, then f\/ix $b_1$ and f\/inally sum over the two permutations of $a_2$, $c_2$.  This results in
		\begin{gather*}
			{\young(\atwo,\ctwo)}
			{\young(\btwo\done\dtwo)}
			g_{ij}
			\Bigl(
				 R\indices{^i_{a_2b_2b_1}}\tilde R\indices{^j_{d_1c_2d_2}}+R\indices{^i_{d_1a_2d_2}}\tilde R\indices{^j_{c_2b_2b_1}} +
				R\indices{^i_{c_2b_2a_2}}\tilde R\indices{^j_{d_1b_1d_2}}+R\indices{^i_{d_1c_2d_2}}\tilde R\indices{^j_{b_1b_2a_2}}\\
\hphantom{{\young(\atwo,\ctwo)}
			{\young(\btwo\done\dtwo)}
			g_{ij}
			\Bigl(}{} +
				 R\indices{^i_{b_1b_2c_2}}\tilde R\indices{^j_{d_1a_2d_2}}+R\indices{^i_{d_1b_1d_2}}\tilde R\indices{^j_{a_2b_2c_2}}
			\Bigr)
			=0 .
		\end{gather*}
	The symmetrisation of this in $b_1$, $b_2$, $d_1$, $d_2$ yields
		\begin{gather*}
			{\young(\atwo,\ctwo)}
			{\young(\bone\btwo\done\dtwo)}
			g_{ij}
			\Bigl(
				 R\indices{^i_{\underline a_2b_2           b_1}}\tilde R\indices{^j_{d_1\underline c_2d_2}}+R\indices{^i_{d_1\underline a_2d_2}}\tilde R\indices{^j_{\underline c_2b_2           b_1}} +
				 R\indices{^i_{\underline c_2b_2\underline a_2}}\tilde R\indices{^j_{d_1           b_1d_2}}+R\indices{^i_{d_1\underline c_2d_2}}\tilde R\indices{^j_{           b_1b_2\underline a_2}}\\
\hphantom{{\young(\atwo,\ctwo)}
			{\young(\bone\btwo\done\dtwo)}
			g_{ij}
			\Bigl(}{}
+
				 R\indices{^i_{           b_1b_2\underline c_2}}\tilde R\indices{^j_{d_1\underline a_2d_2}}+R\indices{^i_{d_1           b_1d_2}}\tilde R\indices{^j_{\underline a_2b_2\underline c_2}}
			\Bigr)
			=0 ,
		\end{gather*}
	where we again underlined antisymmetrised indices.  By the antisymmetry of algebraic curvature tensors in the second index
	pair, all but the fourth and f\/ifth term in the parenthesis vanish and we get
	\[
		{\young(\atwo,\ctwo)}
		{\young(\bone\btwo\done\dtwo)}
		g_{ij}
		\Bigl(
			R\indices{^i_{d_1\underline c_2d_2}}\tilde R\indices{^j_{b_1\underline a_2b_2}}+
			R\indices{^i_{b_1\underline c_2b_2}}\tilde R\indices{^j_{d_1\underline a_2d_2}}
		\Bigr)
		=0.
	\]
	Due to the symmetrisation and antisymmetrisation we can permute the indices of the terms inside the parenthesis to f\/ind that
	\[
		{\young(\atwo,\ctwo)}
		{\young(\bone\btwo\done\dtwo)}
		g_{ij}
		       R\indices{^i_{b_1\underline a_2b_2}}
		\tilde R\indices{^j_{d_1\underline c_2d_2}}
		=0.
	\]
	Recall that this equation has been obtained from \eqref{eq:commutator:hook} by a symmetrisation.  This proves
	\eqref{eq:commutator:hook} $\Rightarrow$ \eqref{eq:commutator:4+2-}.  The converse follows easily by antisymmetrising
	\eqref{eq:commutator:4+2-} in $b_1$, $a_2$, $c_2$.  This achieves the proof of the theorem.
\end{proof}

\begin{Remark}
	Obviously \eqref{eq:commutator:hook} must be true for $\tilde R=R$.  This fact is not evident, but has been proven in~\cite{Schoebel}.
\end{Remark}

For diagonal algebraic curvature tensors $R$ and $\tilde R$ on a four-dimensional vector space the condition~\eqref{eq:commutator:4+2-} can be written in the form
\begin{gather}
	\label{eq:commutator:det}
	\det
	\begin{pmatrix}
		1&R_{ijij}&\tilde R_{ijij}\\
		1&R_{jkjk}&\tilde R_{jkjk}\\
		1&R_{kiki}&\tilde R_{kiki}
	\end{pmatrix}
	=0
\end{gather}
for all distinct $i,j,k\in\{0,1,2,3\}$.

\begin{Corollary}
	\label{cor:Staeckel}
	St\"ackel lines\footnote{Recall Def\/inition~\ref{def:Staeckel-lines}.} in the KS-variety correspond to St\"ackel systems.
	More precisely, the preimage of a St\"ackel line under the map from integrable Killing tensors with diagonal algebraic
	curvature tensor to the KS-variety is a St\"ackel system.  In particular, the Killing tensors in a St\"ackel system have
	simultaneously diagonalisable algebraic curvature tensors.
\end{Corollary}

\begin{proof}
	A St\"ackel line is the projective line through the (collinear but not coinciding) points~$M$, $\iota(\pi(M))$ and
	$\nu(\pi(M))$ in the KS-variety for any non-singular KS-matrix~$M$.  We have to prove that corresponding Killing tensors~$K_0$,~$K_1$ and~$K_2$ mutually commute.  We f\/irst show that~$K_1$ and~$K_2$ commute.  Let $R$ and $\tilde R$ be their
	diagonal algebraic curvature tensors.	From their KS-matrices~\eqref{eq:embeddings} and $\Delta_\alpha=w_\beta-w_\gamma$ we
	can read of\/f the values
	\begin{gather*}
		       w_\alpha =0 ,\qquad      t_\alpha =n_\alpha,\qquad
		\tilde w_\alpha =n_\alpha^2, \qquad \tilde t_\alpha =n_\beta n_\gamma ,
	\end{gather*}
	where we have neglected the scalar curvature, which can be chosen appropriately.  The diagonals of~$R$ and~$\tilde R$ are
	then given by the columns of the matrix
	\[
		\begin{pmatrix}
			1&R_{0101}&\tilde R_{0101}\\
			1&R_{0202}&\tilde R_{0202}\\
			1&R_{0303}&\tilde R_{0303}\\
			1&R_{2323}&\tilde R_{2323}\\
			1&R_{3131}&\tilde R_{3131}\\
			1&R_{1212}&\tilde R_{1212}
		\end{pmatrix}
		=
		\begin{pmatrix}
			1&w_1+t_1&\tilde w_1+\tilde t_1\\
			1&w_2+t_2&\tilde w_2+\tilde t_2\\
			1&w_3+t_3&\tilde w_3+\tilde t_3\\
			1&w_1-t_1&\tilde w_1-\tilde t_1\\
			1&w_2-t_2&\tilde w_2-\tilde t_2\\
			1&w_3-t_3&\tilde w_3-\tilde t_3
		\end{pmatrix}
		=
		\begin{pmatrix}
			1&+n_1&n_1^2+n_2n_3\\
			1&+n_2&n_2^2+n_3n_1\\
			1&+n_3&n_3^2+n_1n_2\\
			1&-n_1&n_1^2-n_2n_3\\
			1&-n_2&n_2^2-n_3n_1\\
			1&-n_3&n_3^2-n_1n_2
		\end{pmatrix}.
	\]
	The four square matrices in \eqref{eq:commutator:det} can be obtained from this matrix by discarding the rows $(1,2,3)$,
	$(1,5,6)$, $(2,4,6)$ respectively $(3,4,5)$ and it is not dif\/f\/icult to see that their determinants vanish.  This shows that
	$K_1$ and $K_2$ commute.

	If $\iota(\pi(M))$ and $\nu(\pi(M))$ are linearly independent, then $M$ is a linear combination of both and hence $K_0$ is a
	linear combination of $g$, $K_1$ and~$K_2$.  In other words, $K_0$, $K_1$ and $K_2$ mutually commute.

	If $\iota(\pi(M))$ and $\nu(\pi(M))$ are linearly dependent, then $n=\pi(M)$ satisf\/ies $\lvert n_1\rvert=\lvert
	n_2\rvert=\lvert n_3\rvert$.  For simplicity, suppose $n_1=n_2=n_3=1$.  The other choices of the signs are analogous.  By
	Lemma~\ref{lem:isokernel}, $M$ is of the form \eqref{eq:isokernel:parametrisation} and thus has $w_\alpha=-t_\alpha$.
	Consequently, the diagonal algebraic curvature tensors $R$ and $\tilde R$ of $K_0$ respectively $K_2$ are given by
	\[
		\begin{pmatrix}
			1&R_{0101}&\tilde R_{0101}\\
			1&R_{0202}&\tilde R_{0202}\\
			1&R_{0303}&\tilde R_{0303}\\
			1&R_{2323}&\tilde R_{2323}\\
			1&R_{3131}&\tilde R_{3131}\\
			1&R_{1212}&\tilde R_{1212}
		\end{pmatrix}
		=
		\begin{pmatrix}
			1&w_1+t_1&\tilde w_1+\tilde t_1\\
			1&w_2+t_2&\tilde w_2+\tilde t_2\\
			1&w_3+t_3&\tilde w_3+\tilde t_3\\
			1&w_1-t_1&\tilde w_1-\tilde t_1\\
			1&w_2-t_2&\tilde w_2-\tilde t_2\\
			1&w_3-t_3&\tilde w_3-\tilde t_3
		\end{pmatrix}
		=
		\begin{pmatrix}
			1&0   &2\\
			1&0   &2\\
			1&0   &2\\
			1&2t_1&0\\
			1&2t_2&0\\
			1&2t_3&0
		\end{pmatrix}.
	\]
	As above, they satisfy the commutation condition \eqref{eq:commutator:det} and we conclude that $K_0$, $K_1$ and $K_2$
	commute.

	It remains to prove that every St\"ackel system consists of Killing tensors with mutually diagonalisable algebraic curvature
	tensors.  This follows from the fact that every St\"ackel system contains a Killing tensor with simple eigenvalues and that
	such a Killing tensor $K$ uniquely determines the St\"ackel system~\cite{Benenti93}.  By the very def\/inition of a~St\"ackel
	system, $K$ is integrable and by Theorem~\ref{thm:slice} we can assume it to have a diagonal algebraic curvature
	tensor.  The corresponding point on the KS variety lies on a~St\"ackel line.  By the above, $K$ then lies in a~St\"ackel system
	with diagonal algebraic curvature tensors.  This proves the statement, since the St\"ackel system determined by $K$ is unique.
\end{proof}

\begin{Remark}
	Our result also shows that every integrable Killing tensor on $S^3$ is contained in a St\"ackel system, not only those with
	simple eigenvalues.  This follows from the above corollary in conjunction with Proposition~\ref{prop:isokernel}.
\end{Remark}

We can now prove Theorem~\ref{thm:blow-up}.  Recall from Proposition~\ref{prop:isokernel} that each St\"ackel
line contains a~unique skew symmetric point $\iota(n)$, where $n\in\mathbb P^2$, and that each non-singular skew symmetric point
$\iota(n)$ determines a unique St\"ackel line, given by the two points $\iota(n)\not=\nu(n)$.  Moreover, the singular skew
symmetric points are the four points $n=(\pm1\!:\!\pm1\!:\!\pm1)\in\mathbb P^2$ for which $\iota(n)=\nu(n)$.  Hence it suf\/f\/ices
to show that the subspace spanned by $\iota(n(t))$ and $\nu(n(t))$ has a well def\/ined limit for $t\to0$ if $n(0)$ is one of
these and that this limit depends on $\dot n(0)\in\mathbb P^2$.  We leave it to the reader to verify that the limit space is
spanned by $\iota(n(0))$ and $\iota(\dot n(0))$.

\section{Geometric constructions of integrable Killing tensors}
\label{sec:geometric}

In this section we present several geometric constructions of integrable Killing tensors and interpret each of them within our
algebraic picture from the last section.

\subsection{Special Killing tensors}

There are two ways in which integrable Killing tensors are related to geodesically equivalent metrics.  The f\/irst is via special
conformal Killing tensors.
\begin{Definition}
	\label{def:SCKT}
	A {\it special conformal Killing tensor} on a Riemannian manifold is a symmetric tensor $L_{\alpha\beta}$ satisfying the
	{\it Sinjukov equation}
	\begin{subequations}
		\label{eq:Sinjukov}
		\begin{gather}
			\nabla_\gamma L_{\alpha\beta}=\lambda_\alpha g_{\beta\gamma}+\lambda_\beta g_{\alpha\gamma},
		\end{gather}
		where
		\begin{gather}
			\lambda=\tfrac12\nabla\tr L
			 ,
		\end{gather}
	\end{subequations}
	as can be seen from contracting $\alpha$ and $\beta$.
\end{Definition}
Special conformal Killing tensors parametrise geodesically equivalent metrics in the following way \cite{Sinjukov}.
\begin{Theorem}
	A metric $\tilde g$ is geodesically equivalent to $g$ if and only if the tensor
	\begin{gather}
		\label{eq:L(g)}
		L\coloneq\left(\frac{\det\tilde g}{\det g}\right)^\frac1{n+1}\tilde g^{-1}
	\end{gather}
	is a special conformal Killing tensor.
\end{Theorem}

From the def\/inition we see that special conformal Killing tensors on a manifold $M$ form a~vector space which is invariant under
the isometry group of $M$.  In other words, they def\/ine a~representation of this group.  The following lemma shows that this
representation is isomorphic to a~subrepresentation of the representation of the isomorphism group on Killing tensors.

\begin{Lemma}
	\label{lem:K(L)}
	If $L$ is a special conformal Killing tensor, then $K\coloneq L-(\tr L)g$ is a Killing tensor.  This defines an injective
	map from the space of special conformal Killing tensors to the space of Killing tensors, which is equivariant with respect
	to the action of the isometry group.
\end{Lemma}
\begin{proof}
	It is straightforward to check the f\/irst statement.  Taking the trace on both sides of $K=L-(\tr L)g$ yields
	\[
		\tr K=(1-n)\tr L
	\]
	and hence
	\begin{gather}
		\label{eq:L(K)}
		L=K-\frac{\tr K}{n-1}g .
	\end{gather}
	This shows injectivity.  Equivariance is obvious.
\end{proof}

\begin{Definition}
	We def\/ine a {\it special Killing tensor} to be a Killing tensor of the form
	\[
		K=L-(\tr L)g,
	\]
	where $L$ is a special conformal Killing tensor.
\end{Definition}

\begin{Proposition}
	A special conformal Killing tensor as well as the corresponding special Killing tensor have vanishing Nijenhuis torsion.  In
	particular, special Killing tensors are integrable.
\end{Proposition}

\begin{proof}
	Substituting the Sinjukov equation \eqref{eq:Sinjukov} into the expression \eqref{eq:Nijenhuis} for $K=L$ shows that the
	Nijenhuis torsion of $L$ is zero.  Substituting $K=L-(\tr L)g$ into \eqref{eq:Nijenhuis} shows that this implies that the
	Nijenhuis torsion of $K$ is also zero.
\end{proof}

For a non-f\/lat constant curvature manifold $M\subset V$ we can identify the subrepresentation of special Killing tensors inside
the space of Killing tensors.  Recall that the space of Killing tensors is isomorphic to the space of algebraic curvature
tensors on $V$ and an irreducible $\GLG(V)$-representation.  Under the subgroup $\SOG(V)\subset\GLG(V)$ this space decomposes
into a Weyl component and a Ricci component.

\begin{Definition}
	\label{def:Kulkarni-Nomizu}
	The {\it Kulkarni--Nomizu product} of two symmetric tensors $h$ and $k$ on a vector space $V$ is the algebraic curvature
	tensor $h\varowedge k$, given by
	\begin{gather}
		\label{eq:Kulkarni-Nomizu}
		(h\varowedge k)_{a_1b_1a_2b_2}
		\coloneq\;
		h_{a_1a_2}k_{b_1b_2}-h_{a_1b_2}k_{b_1a_2}-h_{b_1a_2}k_{a_1b_2}+h_{b_1b_2}k_{a_1a_2} .
	\end{gather}
\end{Definition}

The Ricci component consists of all algebraic curvature tensors of the form $h\varowedge g$, where $h$ is a symmetric tensor.
Therefore the dimension of the Ricci part is equal to $N(N+1)/2$ where $N=\dim V$.  This is exactly the dimension of the space
of special conformal Killing tensors on~$M$~\cite{Sinjukov}.  But there is no other $\SOG(V)$-subrepresentation of the same
dimension inside the space of algebraic curvature tensors.  This shows the following.

\begin{Lemma}
	A special Killing tensor on a non-flat constant curvature manifold is a Killing tensor whose algebraic curvature tensor~$R$
	has a vanishing Weyl tensor.  That is, $R$ is of the form
	\[
		R=h\varowedge g
	\]
	for some symmetric tensor $h$ on $V$, where ``$\,\varowedge\!$'' denotes the Kulkarni--Nomizu product~\eqref{eq:Kulkarni-Nomizu}.
\end{Lemma}

Alternatively, this lemma can be checked directly by verifying that \eqref{eq:L(K)} satisf\/ies \eqref{eq:Sinjukov} if $K$ has an
algebraic curvature tensor of the form $R=h\varowedge g$.

The following consequence is a special case of the so called ``cone construction'' \cite{Matveev&Mounoud}:  Any special
conformal Killing tensor on a Riemannian manifold $M$ can be extended to a covariantly constant symmetric tensor on the metric
cone over $M$.  In our case, where $M\subset V$ is a constant curvature manifold, this cone is nothing but the embedding space
$V$.  That is why a special conformal Killing tensor is the restriction of a constant symmetric tensor on the ambient space.
\begin{Proposition}
	Let $R=h\varowedge g$ be the algebraic curvature tensor of a special Killing tensor $K=L-(\tr L)g$ on a non-flat constant
	curvature manifold $M\subset V$.  Then the corresponding special conformal Killing tensor
	\begin{gather}
		\label{eq:L(K):again}
		L=K-\frac{\tr K}{n-1}g
	\end{gather}
	is the restriction of the $($constant$)$ symmetric tensor
	\begin{gather}
		\label{eq:restriction}
		\hat L=h-\frac{\Tr h}{n-1} g
	\end{gather}
	from $V$ to $M$, where ``\,$\Tr$'' denotes the trace on $V$ and ``\,$\tr$'' the trace on the tangent space.  In particular,
	if $h$ is tracefree, then $L$ is simply the restriction of $h$ from $V$ to $M$.
\end{Proposition}
\begin{proof}
	The special Killing tensor $K$ is given by substituting $R=h\varowedge g$ into \eqref{eq:correspondence} using
	\eqref{eq:Kulkarni-Nomizu}.  This yields
	\[
		K_x(v,w)
		=
		h(x,x)g(v,w)-
		h(x,v)g(x,w)-
		h(x,w)g(x,v)+
		h(v,w)g(x,x)
	\]
	for tangent vectors $v,w\in T_xM$ at a point $x\in M$, i.e.\ vectors $x,v,w\in V$ with $g(x,x)=1$ and $g(x,v)=g(x,w)=0$.
	This simplif\/ies to
	\begin{gather}
		\label{eq:K}
		K_x(v,w)=h(v,w)+h(x,x)g(v,w).
	\end{gather}
	Let $e_1,\ldots,e_n$ be an orthonormal basis of $T_xM$ and complete it with $e_0\coloneq x$ to an orthonormal basis of $V$.
	Then
	\begin{gather}
			\tr K
			 =\sum_{\alpha=1}^nK_x(e_\alpha,e_\alpha)
			 =\sum_{\alpha=1}^n\bigl(h(e_\alpha,e_\alpha)+h(e_0,e_0)g(e_\alpha,e_\alpha)\bigr)\nonumber\\
\hphantom{\tr K}{}
=\sum_{i=0}^nh(e_i,e_i)+(n-1)h(e_0,e_0)
			 =\Tr h+(n-1)h(x,x) .\label{eq:tr(K)}
		\end{gather}
	Recall that ``$\tr$'' denotes the trace on $T_xM$, whereas ``$\Tr$'' denotes the trace on $V$.  Substitu\-ting~\eqref{eq:K}
	and~\eqref{eq:tr(K)} into the right hand side of \eqref{eq:L(K):again} now yields \eqref{eq:restriction}.
\end{proof}

For $S^3$ we can verify the integrability of a special Killing tensor on $S^3$ algebraically, using the Theorems~\ref{thm:slice}
and~\ref{thm:det}:
\begin{Lemma}	\label{lem:special}\quad
	\begin{enumerate} \itemsep=0pt
		\item[$1.$]
			$R=h\varowedge g$ is diagonalisable.  More precisely, $R$ is diagonal in a basis where $h$ is diagonal.  In this
			basis the diagonal elements of $R$ are given in terms of the diagonal elements $h_i\coloneq h_{ii}$ by
			\[
				R_{ijij}=h_i+h_j,
				\qquad
				i\not=j .
			\]
		\item[$2.$]
			If $\dim V=4$, the KS-matrix \eqref{eq:matrix} of $R=h\varowedge g$ is antisymmetric with
			\begin{gather*}
				\Delta_\alpha =0,\qquad
				t_\alpha =h_0+h_\alpha .
			\end{gather*}
			In particular it has determinant zero.
	\end{enumerate}
\end{Lemma}

\begin{Corollary}
	The space of special Killing tensors on $S^3$ corresponds to the projective space of antisymmetric matrices inside the
	KS-variety.
\end{Corollary}

Proposition~\ref{prop:isokernel} and Corollary~\ref{cor:Staeckel} now give a unique representative in each St\"ackel
system.
\begin{Corollary}
	Every St\"ackel system on $S^3$ contains a special Killing tensor which is unique up to multiplication with constants and
	addition of multiples of the metric.
\end{Corollary}

\subsection{Killing tensors of Benenti type}

The second way in which integrable Killing tensors arise from geodesically equivalent metrics is described by the following
theorem \cite{Benenti92, Levi-Civita,Matveev&Topalov98,Painleve}.
\begin{Theorem}
	If a metric $\tilde g$ on a $($pseudo-$)$Riemannian manifold is geodesically equivalent to~$g$, then
	\begin{gather}
		\label{eq:K(g)}
		K\coloneq\left(\frac{\det g}{\det\tilde g}\right)^\frac2{n+1}\tilde g
	\end{gather}
	is an integrable Killing tensor for $g$.
\end{Theorem}

\begin{Definition}
	We will call a Killing tensor of the form \eqref{eq:K(g)} a {\it Benenti--Killing tensor}.
\end{Definition}
Recall that the standard metric on the unit sphere $S^n\subset V$ in a Euclidean vector space $(V,g)$ is the restriction of the
scalar product $g$ from $V$ to $S^n$, which we denoted by $g$ as well.  Its geodesics are the great circles on $S^n$.  Consider
the map
\[
	\begin{array}{@{}rrcll}
		f\colon&S^n&\to&S^n,\\
		&x&\mapsto&f(x)\coloneq\dfrac{Ax}{\lVert Ax\rVert},
	\end{array}
\]
for $A\in\GLG(V)$.  Since $A$ takes hyperplanes to hyperplanes, $f$ takes great circles to great circles.  Hence the pullback
$\tilde g$ of $g$ under $f$ is geodesically equivalent to $g$ by def\/inition.  The Killing tensor obtained from applying
\eqref{eq:K(g)} to $\tilde g$ is given by \cite{Schoebel}
\begin{gather*}
	K_x(v,w) =\frac{g(Ax,Ax)g(Av,Aw)-g(Ax,Av)g(Ax,Aw)}{(\det A)^{\frac4{n+1}}},\qquad
	v,w \in T_xM .
\end{gather*}
Comparing with \eqref{eq:correspondence}, we see that the algebraic curvature tensor of this Killing tensor is proportional to
$h\varowedge h$, where ``$\varowedge$'' denotes the Kulkarni--Nomizu product \eqref{eq:Kulkarni-Nomizu} and
\[
	h(v,w)\coloneq g(Av,Aw)=g\big(v,A^\transpose Aw\big)	 .
\]
Moreover, we can replace the scalar product $h$ by any~-- not necessarily positive def\/inite~-- symmetric tensor.  This motivates
the following def\/inition.
\begin{Definition}
	We will say a Killing tensor on a constant curvature manifold {\it of Benenti type} if its algebraic curvature tensor is
	proportional to
	\[
		R=h\varowedge h
	\]
	for some symmetric tensor $h$ on $V$.
\end{Definition}

A Killing tensor of Benenti type with positive def\/inite $h$ is a Benenti--Killing tensor and thus integrable.  With the aid of
the algebraic integrability conditions one can verify that any Killing tensor of Benenti type is integrable~\cite{Schoebel}.  On
$S^3$ we can check this directly, using the Theorems~\ref{thm:slice} and~\ref{thm:det}:
\begin{Lemma}\label{lem:Benenti}\quad
	\begin{enumerate}\itemsep=0pt
		\item[$1.$]
			For a symmetric tensor $h$ on a vector space $V$, the algebraic curvature tensor $R=h\varowedge h$ is
			diagonalisable.  More precisely, $R$ is diagonal in a basis where $h$ is diagonal.  In this basis the diagonal
			elements of $R$ are given in terms of the diagonal elements $h_i\coloneq h_{ii}$ by
			\begin{gather}
				\label{eq:Benenti:R}
				R_{ijij}=2h_ih_j,
				\qquad
				i\not=j
				 .
			\end{gather}
		\item[$2.$] If $\dim V=4$, the KS-matrix \eqref{eq:matrix} of $R=h\varowedge h$ is given by
			\begin{gather*}
				\Delta_\alpha =(h_0-h_\alpha)(h_\beta-h_\gamma),\qquad
				t_\alpha =h_0h_\alpha-h_\beta h_\gamma
			\end{gather*}
			and has zero determinant and trace.
	\end{enumerate}
\end{Lemma}

\subsection{Benenti systems}

Special Killing tensors and Benenti--Killing tensors only constitute particular cases of a more general construction of
integrable Killing tensors out of special conformal Killing tensors.  This construction appears in the literature in dif\/ferent
guises and under a number of dif\/ferent names:  as ``Newtonian systems of quasi-Lagrangian type''
\cite{Rauch-Wojciechowski&Marciniak&Lundmark}, as ``systems admitting special conformal Killing tensors'' \cite{Crampin03a}, as
``cofactor systems'' \cite{Lundmark}, as ``bi-Hamiltonian structures'' \cite{Blaszak, Ibort&Magri&Marmo}, as
``bi-quasi-Hamiltonian systems'' \cite{Crampin03b, Crampin&Sarlet02}, as ``$L$-systems'' \cite{Benenti05} or as ``Benenti
systems'' \cite{Bolsinov&Matveev}.  Here we follow Bolsinov and Matveev, who call them ``Benenti systems''
\cite{Bolsinov&Matveev}, because their formulation is directly ref\/lected in our algebraic description for $S^3$.

Recall that the adjugate matrix $\Adj M$ of a matrix $M$ is the transpose of its cofactor matrix, which is equal to $(\det
M)M^{-1}$ if $M$ is invertible.  Given a metric, we can identify endomorphisms and bilinear forms.  This allows us to extend
this definition to bilinear forms.  As for endomorphisms, we denote the adjugate of a bilinear form $L$ by $\Adj L$, which is
again a bilinear form.
\begin{Theorem}[\cite{Bolsinov&Matveev}]
	For a special conformal Killing tensor $L$ the family
	\[
		K(\lambda)=\Adj(L-\lambda g),
		\qquad
		\lambda\in{\mathbb R}
	\]
	is a family of mutually commuting integrable Killing tensors.  In particular, for $\lambda=0$ we recover the Benenti--Killing
	tensor~\eqref{eq:K(g)} from~\eqref{eq:L(g)}.
\end{Theorem}
\begin{Definition}
	The family $K(\lambda)=\Adj(L-\lambda g)$ of mutually commuting Killing tensors is called the {\it Benenti system} of the
	special conformal Killing tensor~$L$.
\end{Definition}

In order to state the algebraic counterpart of the above theorem for constant curvature manifolds, we f\/irst need the following
lemma.

\begin{Lemma}
	Let $h$ be a symmetric tensor on a vector space $V$ of dimension $n+1$.  Then the algebraic curvature tensor
	\begin{gather}
		\label{eq:polynomial}
		\frac{(\Adj h)\varowedge(\Adj h)}{\det h}
	\end{gather}
	is well defined and homogeneous in $h$ of degree $n-1$.
\end{Lemma}
\begin{proof}
	Chose a basis in which $h$ is diagonal with diagonal elements $h_i\coloneq h_{ii}$.  Then $\Adj h$ is also diagonal and has
	diagonal elements $h_0\cdots h_{i-1}h_{i+1}\cdots h_n$.  By \eqref{eq:Benenti:R} the algebraic curvature tensor
	\eqref{eq:polynomial} is diagonal and given by
	\begin{gather*}
		R_{ijij}
		 =2
		\frac
			{h_0\cdots h_{i-1}h_{i+1}\cdots h_n\cdot h_0\cdots h_{j-1}h_{j+1}\cdots h_n}
			{h_0\cdots h_n}\\
		\hphantom{R_{ijij}}{}
=2h_0\cdots h_{i-1}h_{i+1}\cdots h_{j-1}h_{j+1}\cdots h_n
	\end{gather*}
	for $i\not=j$. This is well def\/ined and homogeneous in $h$ of degree $n-1$.
\end{proof}

\begin{Proposition}[algebraic representation of Benenti systems]
	\label{prop:R(lambda)}
	Let $L$ be a special conformal Killing tensor on a non-flat constant curvature manifold $M\subset V$, i.e.\ the
	restriction of a~$($constant$)$ symmetric tensor~$\hat L$ from~$V$ to~$M$.  To avoid confusion, we denote the metric on~$V$ by~$\hat g$ and its restriction to $M\subset V$ by~$g$.

	Then the Killing tensors of the Benenti system
	\begin{gather}
		\label{eq:K(lambda)}
		K=\Adj(L-\lambda g),
		\qquad
		\lambda\in{\mathbb R}
	\end{gather}
	correspond to the algebraic curvature tensors
	\begin{gather}
		\label{eq:R(lambda)}
		R=\frac{\Adj(\hat L-\lambda\hat g)\varowedge\Adj(\hat L-\lambda\hat g)}{\det(\hat L-\lambda\hat g)} .
	\end{gather}
	In particular, a member of this Benenti system is of Benenti type unless $\lambda$ is an eigenvalue of~$\hat L$ and it comes
	from a geodesically equivalent metric if~$-\lambda$ is sufficiently large.
\end{Proposition}

\begin{proof}
	For simplicity of notation, we absorb the trace terms ``$\lambda g$'' and ``$\lambda\hat g$'' into~$L$ respectively~$\hat
	L$ by formally setting $\lambda$ to $0$.  At a f\/ixed point $x\in M$ we consider the block decomposition of the endomorphism
	$\mathbf{\hat L}$ under the splitting $V=T_xM\oplus{\mathbb R} x$.	Using the blockwise LU decomposition of~$\mathbf{\hat L}$,
	\[
		\hat{\mathbf L}
		=
		\begin{pmatrix}
			\mathbf L&l\\
			l^\transpose&n
		\end{pmatrix}
		=
		\begin{pmatrix}
			\mathbf L&0\\
			l^\transpose&1
		\end{pmatrix}
		\begin{pmatrix}
			\mathbf \Id&\mathbf L^{-1}l\\
			0&n-l^\transpose\mathbf L^{-1}l
		\end{pmatrix}
		 ,
	\]
	we can compute the block decomposition of the adjugate of $\hat{\mathbf L}$:
	\[
		\Adj\hat{\mathbf L}
		=
		\begin{pmatrix}
			\dfrac1{\det\mathbf L}\Bigl(\det\hat{\mathbf L}\Adj\mathbf L+(\Adj\mathbf L)ll^\transpose(\Adj\mathbf L)^\transpose\Bigr)&-(\Adj\mathbf L)l\\
			-l^\transpose(\Adj\mathbf L)^\transpose&\det\mathbf L
		\end{pmatrix}
		 .
	\]
	Note that we have assumed without loss of generality that~$\mathbf L$ is invertible.  Now take a basis $e_1,\ldots,e_n$ of
	$T_xM$ and complete it with $e_0:=x$ to a basis of~$V$.  With respect to this basis the components of the Killing tensor $K$
	corresponding to the algebraic curvature tensor~\eqref{eq:R(lambda)} read
	\[
		K_{ij}
		=R_{0i0j}
		=\frac{(\Adj\hat L)_{00}(\Adj\hat L)_{ij}-(\Adj\hat L)_{0i}(\Adj\hat L)_{0j}}{\det\hat L}
		=(\Adj L)_{ij} .
	\]
	That is, $K$ is given by \eqref{eq:K(lambda)}, as was to be shown.
\end{proof}

The adjugate of an $(n\times n)$-matrix $M$ is a polynomial in $M$ of degree $n-1$, given by
\[
	\Adj M=\left.\frac{\chi_M(\lambda)-\chi_M(0)}\lambda\right\rvert_{\lambda=M} ,
\]
where $\chi_M$ is the characteristic polynomial of $M$.  Therefore a Benenti system is a polynomial in~$\lambda$ of degree~$n-1$:
\[
	K(\lambda)
	=\Adj(L-\lambda g)
	=\sum_{j=0}^{n-1}K_j\lambda^j .
\]
Up to signs, the leading coef\/f\/icient of this polynomial is the metric, the next coef\/f\/icient is the special Killing tensor~\eqref{eq:L(g)} and the constant term is the Benenti--Killing tensor~\eqref{eq:K(g)}:
\begin{gather*}
	K_{n-1}
	=\lim_{\lambda\to\infty}\frac1{\lambda^{n-1}}K(\lambda)
	=\lim_{\lambda\to\infty}\Adj\bigl(\tfrac1\lambda L-g\bigr)
	=(-1)^{n-1}g,\\
	K_{n-2}
	=\left.\frac d{d\lambda}\right|_{\lambda=0}\lambda^{n-1}K\left(\tfrac 1\lambda\right)
	=\left.\frac d{d\lambda}\right|_{\lambda=0}\Adj\left(\lambda L-g\right)
	=(-1)^{n-1}\bigl(L-(\tr L)g\bigr),\\
	K_0 =K(0)=\Adj L.
\end{gather*}
Note that since all the $K(\lambda)$ commute, the $n$ coef\/f\/icients~$K_j$ span a vector space of mutually commuting integrable
Killing tensors.  If the $n$ coef\/f\/icients are linearly independent, then this space is a St\"ackel system and the $K_j$ form a~basis therein.  However, it is in general an open question under which conditions this will be the case.  For $S^3$ we can
derive the answer directly from our algebraic description.

\begin{Lemma}
	\label{lem:Adjstar}
	For a symmetric tensor $h$ on a four dimensional vector space we have
	\begin{gather}
		\label{eq:Adjstar}
		\frac{(\Adj h)\varowedge(\Adj h)}{\det h}=\Hodge(h\varowedge h)\Hodge
		 ,
	\end{gather}
	where $\Hodge$ is the Hodge star operator.
\end{Lemma}
\begin{Remark}
	Regarding \eqref{eq:branching}, the conjugation of an algebraic curvature tensor with the Hodge star is simply a sign change
	in the trace free Ricci part:
	\begin{gather*}
		\mathbf R =\mathbf W+\mathbf T+\mathbf S, \qquad
		*\mathbf R* =\mathbf W-\mathbf T+\mathbf S.
	\end{gather*}
\end{Remark}
\begin{proof}
	Choose a basis in which $h$ is diagonal with diagonal elements $h_i\coloneq h_{ii}$.  For simplicity of notation we assume
	that $\det h\not=0$.  Then $h_i\not=0$ and $\Adj h$ is diagonal with diagonal elements $\det h/h_i$.  Substituting $\Adj h$
	for $h$ in \eqref{eq:Benenti:R} and comparing to \eqref{eq:Benenti:R} shows that the diagonal algebraic	curvature tensors
	\begin{gather*}
		\tilde R =(\Adj h)\varowedge(\Adj h),\qquad
		R =h\varowedge h
	\end{gather*}
	are related by
	\[
		\tilde R_{ijij}
		=2\frac{\det h}{h_i}\frac{\det h}{h_j}
		=2(\det h)\frac{h_0h_1h_2h_3}{h_ih_j}
		=(\det h)2h_kh_l
		=(\det h)R_{klkl}
		 ,
	\]
	for $\{i,j,k,l\}=\{0,1,2,3\}$.  Now the identity~\eqref{eq:Adjstar} follows from~\eqref{eq:chiral} and the fact that $R$ and
	$\Hodge R\Hodge$ dif\/fer in the sign of the trace free Ricci part.
\end{proof}

\begin{Proposition}
	\label{prop:triple}
	Let $L$ be a special conformal Killing tensor on $S^3\subset V$, i.e.\ the restriction of a~$($constant$)$ symmetric tensor
	$\hat L$ from $V$ to $S^3$, and consider the corresponding Benenti system
	\[
		K(\lambda)=\Adj(L-\lambda g)=K_2\lambda^2+K_1\lambda+K_0.
	\]
	Then the coefficients $K_0$, $K_1$ and $K_2$ span a St\"ackel system if and only if $\hat L$ has no triple eigenvalue.
\end{Proposition}

\begin{proof}
	By Proposition~\ref{prop:R(lambda)} and the identity \eqref{eq:Adjstar} the algebraic curvature tensor of $K(\lambda)$
	is the polynomial
	\begin{gather*}
		R(\lambda)
		 =\Hodge\bigl((\hat L-\lambda g)\varowedge(\hat L-\lambda g)\bigr)\Hodge
		 =\lambda^2 g\varowedge g+2\lambda\hat L\varowedge g+\Hodge\big(\hat L\varowedge\hat L\big)\Hodge
		 \eqcolon\lambda^2R_2+\lambda R_1+R_0		 .
	\end{gather*}
	The coef\/f\/icients $R_0$, $R_1$ and $R_2$ are the algebraic curvature tensors of $K_0$, $K_1$ and $K_2$ respectively.  We can
	assume them to be diagonal.  Since $K_2$ is minus the metric, $K_0$, $K_1$ and $K_2$ are linearly independent if and only if
	$R_0$ and $R_1$ def\/ine two distinct points in the KS-variety.

	Choose a basis in which $\hat L$ is diagonal with the eigenvalues $\Lambda_i$ on the diagonal.  Without loss of generality
	we can assume that $\hat L$ is trace free.  Motivated by \eqref{eq:iso:diagonal}, we can thus parametrise~$\hat L$ by a
	vector $n\in{\mathbb R}^3$ as
	\begin{subequations}
		\begin{gather}
			\label{eq:Lambda(n)}
			\Lambda_0      =\frac{n_\alpha+n_\beta+n_\gamma}2, \qquad
			\Lambda_\alpha =\frac{n_\alpha-n_\beta-n_\gamma}2
			 ,
		\end{gather}
		where
		\begin{gather}
			\label{eq:n(Lambda)}
			n_\alpha=\Lambda_0+\Lambda_\alpha=-(\Lambda_\beta+\Lambda_\gamma) .
		\end{gather}
		\end{subequations}
	Then $\hat L$ has a triple eigenvalue if and only if $\lvert n_1\rvert=\lvert n_2\rvert=\lvert n_3\rvert$.  This condition
	is equivalent to the condition that the KS-matrices $\nu(n)$ and $\iota(n)$ in~\eqref{eq:embeddings} are proportional,
	i.e.\ def\/ine the same point in the KS-variety.  Our proof will therefore be f\/inished if we show that $R_0$ and $R_1$
	def\/ine the points $\nu(n)$ respectively $\iota(n)$ in the KS-variety.

	That the algebraic curvature tensor $R_1=2\hat L\varowedge g$ def\/ines the point $\iota(n)$ is the statement of
	Lemma~\ref{lem:special}.  And since $R_0=\Hodge(\hat L\varowedge\hat L)\Hodge$ dif\/fers from $\hat L\varowedge\hat L$ in
	the sign of the trace free Ricci part, Lemma~\ref{lem:Benenti} implies that for $R_0$ the KS-matrix \eqref{eq:matrix}
	is given by
	\begin{gather*}
		\Delta_\alpha =(\Lambda_0-\Lambda_\alpha)(\Lambda_\beta-\Lambda_\gamma) =n_\beta^2-n_\gamma^2,\qquad
		t_\alpha =-(\Lambda_0\Lambda_\alpha-\Lambda_\beta\Lambda_\gamma) =n_\beta n_\gamma .
	\end{gather*}
	Now compare this to $\nu(n)$ in \eqref{eq:embeddings} to f\/inish the proof.
\end{proof}

The Benenti system of a special conformal Killing tensors $L$ def\/ines a map ${\mathbb R}\to\mathcal K(S^3)$ given by $\lambda\mapsto
K(\lambda)=\Adj(L-\lambda g)$.  In the proof above we have actually shown that this def\/ines a map from $\mathbb
P^1={\mathbb R}\cup\{\infty\}$ to the KS-variety which maps $0\mapsto\nu(n)$ and $\infty\mapsto\iota(n)$, where $n\in\mathbb P^2$ is
determined by $L$.  If the points $\iota(n)$ and $\nu(n)$ coincide, the image of this map is a single point.  If they are
distinct, the image is -- by def\/inition -- the St\"ackel line through both points.  In this case the image of the Benenti system
is this St\"ackel line without the point $\iota(n)$.

We can now prove Theorem~\ref{thm:parametrisation}.  It suf\/f\/ices to show that the image of diagonal algebraic curvature
tensors of the form $\lambda_2h\varowedge h$ is dense in the KS-variety.  By Proposition~\ref{prop:R(lambda)} the algebraic
curvature tensors of this form are dense in a Benenti system, so it suf\/f\/ices to show that the image of Benenti systems in the
KS-variety is dense.  But this follows from the preceding paragraph and part~2 
of
Proposition~\ref{prop:isokernel}.

\subsection{Extension of integrable Killing tensors}
\label{sec:extension}

If we f\/ix an orthogonal decomposition ${\mathbb R}^{n_1+n_2+2}={\mathbb R}^{n_1+1}\oplus{\mathbb R}^{n_2+1}$, then two algebraic curvature tensors $R_1$ on
${\mathbb R}^{n_1+1}$ and $R_2$ on ${\mathbb R}^{n_2+1}$ def\/ine an algebraic curvature tensor $R_1\oplus R_2$ on ${\mathbb R}^{n_1+n_2+2}$.
Correspondingly, two Killing tensors $K_1$ on $S^{n_1}$ and $K_2$ on $S^{n_2}$ def\/ine a Killing tensor $K_1\oplus K_2$ on
$S^{n_1+n_2+1}$.  In particular, we can extend any Killing tensor $K$ on $S^{n-1}$ by zero to a Killing tensor $K\oplus0$ on
$S^n$.

From the algebraic integrability conditions we see that $K_1\oplus K_2$ is integrable if $K_1$ and $K_2$ are both integrable.
This direct sum operation def\/ines an embedding
\[
	\mathcal K(S^{n_1})\times\mathcal K(S^{n_2})\hookrightarrow\mathcal K\big(S^{n_1+n_2+1}\big),
\]
where $\mathcal K(S^n)$ denotes the space of integrable Killing tensors on $S^n$.  In particular we have an embedding
\[
	\mathcal K\big(S^{n-1}\big)\hookrightarrow\mathcal K\big(S^n\big).
\]

\begin{Definition}
	We will call (integrable) Killing tensors arising in the way described above {\it extensions} of Killing tensors from
	$(S^{n_1},S^{n_2})$ respectively $S^{n-1}$.
\end{Definition}

As an example, let us determine which integrable Killing tensors on $S^3$ are extensions from~$S^2$ respectively~$(S^1, S^1)$.
Note that in dimension one and two all Killing tensors are integrable.  This follows directly from the integrability conditions.

An algebraic curvature tensor $R$ on ${\mathbb R}^4$ is an extension by zero of an algebraic curvature tensor on ${\mathbb R}^3$ under the
inclusion $(x_1,x_2,x_3)\mapsto(0,x_1,x_2,x_3)$ if and only if $R_{a_1b_1a_2b_2}$ is zero whenever one of the indices is zero.
In particular, the diagonal algebraic curvature tensor \eqref{eq:diagonal} on $S^3$ is an extension from $S^2$ if and only if
$w_\alpha+\frac s{12}=-t_\alpha$.  In this case the KS-matrix \eqref{eq:matrix} is equal to \eqref{eq:isokernel:parametrisation}
with $n_1=n_2=n_3$.  Therefore the extensions of Killing tensors from $S^2$ correspond to the isokernel planes in the
KS-variety.  In particular, the face centers of the octahedron in the KS-variety correspond to the extension of the metric on
$S^2$ under the four embeddings $S^2\subset S^3$ given by intersecting $S^3\subset{\mathbb R}^4$ with each of the four coordinate
hyperplanes.

An algebraic curvature tensor $R$ on ${\mathbb R}^4$ is a sum under the decomposition of $(x_0,x_1,x_2,x_3)$ into $(x_0,x_1)$ and
$(x_2,x_3)$ if and only if the components $R_{a_1b_1a_2b_2}$ are zero unless all four indices are either in $\{0,1\}$ or in
$\{2,3\}$.  In particular, a diagonal algebraic curvature tensor \eqref{eq:diagonal} is an extension from $(S^1,S^1)$ if and
only if $w_2+\frac s{12}=t_2=w_3+\frac s{12}=t_3=0$.  In this case the KS-matrix \eqref{eq:matrix} is equal to
\[
	\begin{pmatrix}
		0&0&0\\
		0&-w_1&-t_1\\
		0&+t_1&+w_1\\
	\end{pmatrix}
	=
	-\frac{w_1-t_1}2V_{-1}
	-\frac{w_1+t_1}2V_{+1}
	 ,
\]
where $V_{+1}$ and $V_{-1}$ are rank one singular points in the KS-variety, cf.\ \eqref{eq:singular}.  Therefore the
extensions of Killing tensors from $(S^1,S^1)$ correspond to the projective lines joining opposite vertices of the octahedron in
the KS-variety.  In particular, each vertex of the octahedron corresponds to an extension of the metrics on two orthogonal
copies of $S^1$ in $S^3$, given by intersecting $S^3\subset{\mathbb R}^4$ with a pair of orthogonal coordinate planes.

\section{Separation coordinates}
\label{sec:coordinates}

We eventually demonstrate how separation coordinates on $S^3$ and their classif\/ication arise naturally and in a purely algebraic
way from our description of integrable Killing tensors.

Kalnins and Miller proved that separation coordinates on $S^n$ are always orthogonal \cite{Kalnins&Miller86}.  From the work of
Eisenhart we know that orthogonal separation coordinates are in bijective correspondence with St\"ackel systems \cite{Eisenhart}.
And in the preceding sections we have proven that on $S^3$ every St\"ackel system contains a special Killing tensor which is
essentially unique.  By def\/inition, special Killing tensors are in bijective correspondence with special conformal Killing
tensors.  But Crampin has shown that the eigenvalues of a special conformal Killing tensor are constant on the coordinate
hypersurfaces of the corresponding separation coordinates \cite{Crampin03a}.  Hence, if these eigenvalues are simple and
non-constant, they can be used as separation coordinates.  This reduces the classif\/ication of separation coordinates on $S^3$ to
a computation of the eigenvalues of special conformal Killing tensors.  The computation is considerably simplif\/ied by the fact
that a special conformal Killing tensor on $S^n\subset V$ is the restriction of a (constant) symmetric tensor on~$V$.

\subsection{Eigenvalues of special conformal Killing tensors}

Let $L$ be a special conformal Killing tensor on a non-f\/lat constant curvature manifold $M\subset V$, i.e.\ the restriction
from $V$ to $M$ of a (constant) symmetric tensor $\hat L$.  From now on we consider $L$ and $\hat L$ as endomorphisms.  We want
to compute the eigenvalues $\lambda_1(x)\le\cdots\le\lambda_n(x)$ of $L\colon T_xM\to T_xM$ for a f\/ixed point $x\in M$ and
relate them to the (constant) eigenvalues $\Lambda_0\le\cdots\le\Lambda_n$ of $\hat L\colon V\to V$.

\begin{figure}[h]
	\centering
	\includegraphics[width=140mm]{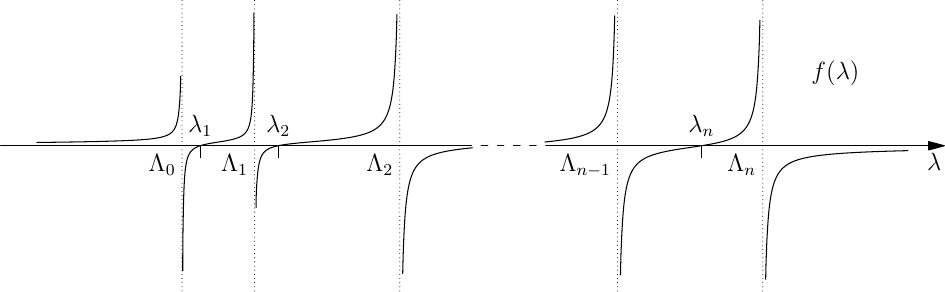}
	\caption{Eigenvalues of a special conformal Killing tensor on $S^n$.}
	\label{fig:eigenvalues}
\end{figure}

Let $P:V\to V$ be the orthogonal projection from $V\isom T_xM\oplus{\mathbb R} x$ to $T_xM$.  Then the restriction of $P\hat LP$ to
$T_xM$ is $L$ and the restriction of $P\hat LP$ to $x$ is zero.  Therefore $x$ is an eigenvector of $P\hat LP$ with eigenvalue
$0$ and the remaining eigenvectors and eigenvalues are those of $L$.  This means the eigenvalue equation $Lv=\lambda v$ is
equivalent to $P\hat LPv=\lambda v$ for $v\perp x$.  By the def\/inition of $P$ we have $Pv=v$ and $P\hat Lv=\hat Lv-g(\hat
Lv,x)x$.  Hence we seek common solutions to the two equations
\begin{gather*}
	(\hat L-\lambda)v =g(\hat Lv,x)x,\qquad g(x,v)=0.
\end{gather*}
If $\lambda$ is not an eigenvalue of $\hat L$, then $(\hat L-\lambda)$ is invertible and $v$ is proportional to $(\hat
L-\lambda)^{-1}x$.  The condition $v\perp x$ then yields the equation
\[
	g\bigl(\big(\hat L-\lambda\big)^{-1}x,x\bigr)=0
\]
for $\lambda$.  In an eigenbasis of $\hat L$ this equation reads
\begin{gather}
	\label{eq:eigenvalues}
	q(\lambda)\coloneq\sum_{k=0}^n\frac{x_k^2}{\Lambda_k-\lambda}=0 .
\end{gather}
That is, the zeroes of the function $q(\lambda)$ are eigenvalues of~$L$.  As depicted in Fig.~\ref{fig:eigenvalues},
$q(\lambda)$ goes to~$0$ for~$\lambda\to\pm\infty$, has poles at the eigenvalues~$\Lambda_k$ of $\hat L$ unless $x_k=0$ and is
monotonely increasing in between.

Let us f\/irst suppose that $\hat L$ has only simple eigenvalues.  In this case \eqref{eq:eigenvalues} is the def\/ining equation
for elliptic coordinates $\lambda_1(x),\ldots,\lambda_n(x)$ on $S^n$ with parameters $\Lambda_0<\Lambda_1<\cdots<\Lambda_n$~\cite{Neumann}.  If $x_k\not=0$ for all $k=0,\ldots,n$, then $q$ has exactly $n+1$ poles at $\Lambda_0,\ldots,\Lambda_n$ and $n$
zeroes in between.  If $x_k=0$, then $q$ has no pole at $\Lambda_k$.  But $x_k=0$ means that the $k$-th eigenvector of~$\hat L$
is orthogonal to $x$ and hence also an eigenvector of $L$ with eigenvalue~$\Lambda_k$.  We conclude that the eigenvectors~$v_k$
of $\hat L$ which are orthogonal to $x$ are eigenvectors of $L$ with the same eigenvalue $\Lambda_k$, and that the remaining
eigenvalues of $L$ are the zeroes $\lambda_k$ of the function $q(\lambda)$ with corresponding eigenvector $v_k=(\hat
L-\lambda_k)^{-1}x$.  In particular, the eigenvalues of~$\hat L$ and~$L$ intertwine:
\begin{gather}
	\label{eq:squeezing}
	\Lambda_0\le\lambda_1(x)\le\Lambda_1\le\lambda_2(x)\le\cdots\le\lambda_n(x)\le\Lambda_n
	 .
\end{gather}
Note that the common denominator of the fractions in the sum \eqref{eq:eigenvalues} is the characteristic polynomial $\chi_{\hat
L}(\lambda)$ of $\hat L$.  Therefore the function $q$ is the quotient of the characteristic polynomials of~$L$ and~$\hat L$.

If $\hat L$ has a simple eigenvalue, then the corresponding eigenspace is one dimensional and has only a trivial intersection
with a generic tangent space to $M$, because $M\subset V$ is of codimension one.  If $\hat L$ has a double eigenvalue, then the
corresponding eigenspace is a plane and has a~one-dimensional intersection with a generic tangent space to $M$.  This def\/ines a
f\/ield of common eigenvectors for~$L$ and $\hat L$ on $M$.  The corresponding eigenvalue is the same for $L$ and $\hat L$ and
hence constant on~$M$.  This also follows from \eqref{eq:squeezing}.  In the same way an eigenvalue of $\hat L$ with
multiplicity~\mbox{$m+1$} def\/ines an $m$-dimensional eigenspace distribution of $L$ with constant eigenvalue.
\begin{Proposition}
	\label{prop:eigenvalues}
	Let $L$ be a special conformal Killing tensor on a constant curvature manifold $M\subset V$, i.e.\ the restriction of a
	constant symmetric trace free tensor $\hat L$ on $V$.  Then:
	\begin{enumerate}\itemsep=0pt
		\item[$1.$]
			If $\hat L$ has only simple eigenvalues, then the eigenvalues of $L$ define elliptic coordinates on $M$ whose
			parameters are the eigenvalues of $\hat L$.
		\item[$2.$] Any multiple eigenvalue of $\hat L$ is a constant eigenvalue for $L$.
		\item[$3.$] $L$ has simple eigenvalues if and only if $\hat L$ has at most double eigenvalues.
	\end{enumerate}
\end{Proposition}

\subsection[Killing tensors and separation coordinates on $S^2$]{Killing tensors and separation coordinates on $\boldsymbol{S^2}$}
\label{sec:S2}

Before we turn to $S^3$, let us brief\/ly discuss integrable Killing tensors and separation coordinates on~$S^2$.  The reason is
that their extensions will appear on~$S^3$.  First note that in dimension two every Killing tensor is integrable.  So we do not
have to deal with the integrability conditions.  Note also that algebraic curvature tensors on~${\mathbb R}^3$ have a zero Weyl part and
are thus determined by their Ricci tensor~$T$ alone.  This means that Killing tensors on $S^2\subset{\mathbb R}^3$ correspond to
symmetric tensors~$T$ on~${\mathbb R}^3$, where the metric corresponds to the identity.

St\"ackel systems on $S^2$ are of dimension two and contain the metric.  That implies that each St\"ackel system on $S^2$
corresponds to a plane of Ricci tensors spanned by the identity and some trace free symmetric tensor, the latter being unique up
to multiples.  Moreover, the isometry group acts by conjugation on the Ricci tensor.  Accordingly, the classif\/ication of
separation coordinates on $S^2$ modulo isometries is equivalent to the classif\/ication of trace free symmetric tensors $T$ on
${\mathbb R}^3$ under the orthogonal group $\OG(3)$.

In particular, every Killing tensor on $S^2$ is special and we can make use of the results in the previous section.  Depending
on the multiplicities of the eigenvalues $t_1\le t_2\le t_3$ of the trace free Ricci tensor $T$ of a Killing tensor on $S^2$ we
consider the following two cases.

{\bf Elliptic coordinates.}
If $t_1<t_2<t_3$, then the restriction $L$ of $T$ from ${\mathbb R}^3$ to $S^2\subset{\mathbb R}^3$ has two dif\/ferent eigenvalues, given by the
zeros of its characteristic polynomial
\[
	\chi_L(\lambda)=
	(t_2-\lambda)(t_3-\lambda)x_1^2+
	(t_3-\lambda)(t_1-\lambda)x_2^2+
	(t_1-\lambda)(t_2-\lambda)x_3^2=0	 .
\]
Using that $x_1^2+x_2^2+x_3^2=1$ and $t_1+t_2+t_3=0$ and writing the endomorphisms $T$ and $\Adj T$ as forms, this quadratic
equation reads
\[
	\lambda^2+T(x,x)\lambda+(\Adj T)(x,x)=0	 .
\]
The two solutions $\lambda_1(x)$ and $\lambda_2(x)$ are real and satisfy $t_1<\lambda_1(x)<t_2<\lambda_2(x)<t_3$.  As mentioned
above, they def\/ine elliptic coordinates on $S^2$ with parameters $t_1<t_2<t_3$.

{\bf Spherical coordinates.}
If $t_1<t_2=t_3$, the characteristic polynomial of $L$ reduces to
\[
	\chi_L(\lambda)=
	(t_2-\lambda)
	\bigl[
		(t_2-\lambda)x_1^2+
		(t_1-\lambda)\big(x_2^2+x_3^2\big)
	\bigr]
	=0
	 .
\]
One of the eigenvalues only depends on $x_1$ and the other is constant:
\begin{gather*}
	\lambda_1(x) =t_2x_1^2+t_1\big(1-x_1^2\big),\qquad
	\lambda_2(x) =t_2.
\end{gather*}
The curves of constant $\lambda_1(x)$ are circles of constant latitude and the meridians complete them to an orthogonal
coordinate system~-- the spherical coordinates.

The case $t_1=t_2<t_3$ is analogue and yields spherical coordinates as well.  Note that even so, both coordinate systems are not
isometric.  This is because they dif\/fer not only by an isometry but also by an exchange of the two coordinates.

We can use the action of the isometry group $\OG(3)$ to diagonalise the Ricci tensor and thereby the algebraic curvature tensor.
Hence the projective variety of integrable Killing tensors on~$S^2$ with diagonal algebraic curvature tensor is isomorphic to
the projective plane~$\mathbb P^2$, where the point $(1:1:1)\in\mathbb P^2$ corresponds to the metric on~$S^2$.  Moreover,
St\"ackel systems on~$S^2$ correspond to projective lines through this point, since they always contain the metric.  The residual
isometry group action on Killing tensors with diagonal algebraic curvature tensor is the natural action of~$S_3$ on~$\mathbb
P^2$ by permutations.

\subsection{The classif\/ication}

We now recover the well known classif\/ication of separation coordinates on $S^3$ \cite{Eisenhart,Kalnins&Miller86}.  In
particular, this will explain the graphical procedure given in \cite{Kalnins&Miller86} in terms of the classif\/ication of trace
free symmetric tensors under the orthogonal group.  The results are summarised in Table~\ref{tab:classification}.

\begin{table}[t]	\centering
	\caption{Classif\/ication of separation coordinates on $S^3$.}
	\label{tab:classification}
\vspace{1mm}

	\begin{tabular}{lllcc}
		\toprule
		type       & coordinates             & induced from    & \cite{Eisenhart} & \cite{Kalnins&Miller86} \\
		\midrule
		$(0123)    $  & elliptic                & --               &   V & (1) \\
		$(01(23))  $  & oblate Lam\'e rotational  & --              &   I & (2a) \\
		$(0(12)3)  $  & prolate Lam\'e rotational & --               &   I & (2b) \\
		$((01)(23))$  & cylindrical             & $(S^1,S^1)$     & III & (5) \\
		$(0(123))  $  & Lam\'e subgroup reduction & $S^2$ elliptic  &  IV & (3) \\
		$(0(1(23)))$  & spherical               & $S^2$ spherical &  II & (4) \\
		\bottomrule
	\end{tabular}

\end{table}

As before, we will denote by $\Lambda_0\le\Lambda_1\le\Lambda_2\le\Lambda_3$ the ordered eigenvalues of the symmetric tensor
$\hat L$ associated to a St\"ackel system and parametrised by a vector $(n_1,n_2,n_3)\in{\mathbb R}^3$ via~\eqref{eq:Lambda(n)}.  Recall
that we have the freedom to add a multiple of the metric to a Killing tensor or to multiply it with a constant.  Both operations
only af\/fect the eigenvalues, not the eigenspaces, and therefore result in dif\/ferent parametrisations of the same coordinate
system.  Whereas Kalnins and Miller use this freedom to set $\Lambda_0=0$ and $\Lambda_1=1$, here we use it to assume that $\hat
L$ is trace free and then to consider $(\Lambda_0:\Lambda_1:\Lambda_2:\Lambda_3)$ respectively $(n_1:n_2:n_3)$ as a point in~$\mathbb P^2$.

\looseness=-1
In particular, changing the sign of $\hat L$ reverses the order of the eigenvalues.  This leaves us with the following f\/ive
dif\/ferent alternatives for the multiplicities of the ordered set of eigenvalues of~$\hat L$.  For simplicity we write~``$i$''
for~``$\Lambda_i$'' and denote multiple eigenvalues by parenthesising them:
\begin{gather*}
	(0123),\qquad
	(01(23)),\qquad
	(0(12)3),\qquad
	((01)(23)),\qquad	(0(123)) .
\end{gather*}

\looseness=-1
{\bf (0123) Elliptic coordinates.}
We have already seen that if the eigenvalues of~$\hat L$ are simple, then the eigenvalues of $L$ are also simple and def\/ine
elliptic coordinates on~$S^n$ given implicitly by~\eqref{eq:eigenvalues}.

{\bf (01(23)) Oblate Lam\'e rotational coordinates.}
If $\hat L$ has a double eigenvalue $\Lambda_{k-1}=\Lambda_k$, then the corresponding eigenvalue $\lambda_k$ of $L$ is constant
and can not be used as a coordinate.  Nevertheless, we can derive the missing coordinate from the fact that the coordinate
hypersurfaces are orthogonal to the eigenspaces of $L$.  Let $E$ be an eigenplane of $\hat L$ and denote by $E^\perp$ its
orthogonal complement in $V$.  Recall that the corresponding eigenspace of $L$ at a point $x\in S^n$ is given by the
intersection of $E$ with the tangent space $T_xS^n$.  Its orthogonal complement in $T_xS^n$ is the intersection of $T_xS^n$ with
the hyperplane that contains both $x$ and $E^\perp$.  The coordinate hypersurfaces are therefore the intersections of
hyperplanes containing $E^\perp$ with $S^n$ and can be parametrised by the polar angle in $E$.  Hence this angle can be taken as
a coordinate.

\looseness=-1
By \eqref{eq:squeezing}, $L$ has single eigenvalues given by \eqref{eq:eigenvalues} if and only if $\hat L$ has at most double
eigenvalues and every double eigenvalue of $\hat L$ determines a constant eigenvalue for $L$.  If we replace every constant
eigenvalue by the polar angle in the corresponding eigenplane of $\hat L$, we obtain a complete coordinate system.  Concretely,
in the case $\Lambda_0<\Lambda_1<\Lambda_2=\Lambda_3$ the separation coordinates are def\/ined by the solutions $\lambda_1(x)$ and
$\lambda_2(x)$ of~\eqref{eq:eigenvalues} together with the polar angle in the $(x_2,x_3)$-plane.

For $\Lambda_0<\Lambda_1<\Lambda_2=\Lambda_3$ we have $n_1<n_2=n_3$.  The corresponding St\"ackel line thus goes through the two
points in the KS-variety given by the KS-matrices \eqref{eq:embeddings}:
\begin{gather*}
	\begin{pmatrix}
		0&-n_2&n_2\\
		n_2&0&-n_1\\
		-n_2&n_1&0
	\end{pmatrix},\qquad
	\begin{pmatrix}
		0&-n_1n_2&n_1n_2\\
		n_1n_2&n_2^2-n_1^2&-n_2^2\\
		-n_1n_2&n_2^2&n_1^2-n_2^2
	\end{pmatrix}	 .
\end{gather*}
Subtracting $n_1$ times the f\/irst matrix from the second, we see that the same St\"ackel line is also determined by the two points
\begin{gather*}
	\begin{pmatrix}
		0&-n_2&n_2\\
		n_2&0&-n_1\\
		-n_2&n_1&0
	\end{pmatrix},\qquad
	V_{-1} =
	\begin{pmatrix}
		0&0&0\\
		0&+1&-1\\
		0&+1&-1
	\end{pmatrix}
	 .
\end{gather*}

{\bf (0(12)3) Prolate Lam\'e rotational coordinates.}
The case $\Lambda_0<\Lambda_1=\Lambda_2<\Lambda_3$ is analogue and yields Lam\'e rotational coordinates which dif\/fer from the
above by an isometry and an exchange of two of the coordinates.

{\bf ((01)(23)) Cylindrical coordinates.}
In the case $\Lambda_0=\Lambda_1<\Lambda_2=\Lambda_3$ equation \eqref{eq:eigenvalues} has only one non-constant solution, namely
$\lambda_2(x)$:
\begin{gather*}
	\lambda_1 =\Lambda_1 \le
	\lambda_2(x) =\Lambda_1\big(x_2^2+x_3^2\big)+\Lambda_3\big(x_0^2+x_1^2\big) \le
	\lambda_3 =\Lambda_3.
\end{gather*}
The separation coordinates are obtained by replacing the constant eigenvalues $\lambda_1$ and $\lambda_3$ by the polar angle in
the $(x_0,x_1)$ respectively $(x_2,x_3)$-plane.

For $\Lambda_0=\Lambda_1<\Lambda_2=\Lambda_3$ we have $n_1<n_2=n_3=0$ and the corresponding St\"ackel line goes through the two
points in the KS-variety given by the KS-matrices \eqref{eq:embeddings}:
\begin{gather*}
	\begin{pmatrix}
		0&0&0\\
		0&0&-n_1\\
		0&+n_1&0
	\end{pmatrix},\qquad
	\begin{pmatrix}
		0&0&0\\
		0&-n_1^2&0\\
		0&0&n_1^2
	\end{pmatrix} .
\end{gather*}
Equivalently, this St\"ackel line goes through the two points
\begin{gather*}
	V_{-1} =
	\begin{pmatrix}
		0&0&0\\
		0&+1&-1\\
		0&+1&-1
	\end{pmatrix},\qquad
	V_{+1} =
	\begin{pmatrix}
		0&0&0\\
		0&+1&+1\\
		0&-1&-1
	\end{pmatrix}
	 ,
\end{gather*}
i.e.\ opposed vertices of the octahedron in the KS-variety.  We have seen that this describes extensions from~$(S^1,S^1)$ to~$S^3$.

\looseness=-1
{\bf (0(123)) Coordinates extended from $S^2\subset S^3$.}
In the remaining case $\Lambda_0\!<\!\Lambda_1\!=\!\Lambda_2\!=\!\Lambda_3$, equation \eqref{eq:eigenvalues} has the solutions
\begin{gather*}
	\lambda_1(x) =\Lambda_1x_0^2+\Lambda_0\big(x_1^2+x_2^2+x_3^2\big),\qquad
	\lambda_2 =\lambda_3=\Lambda_1.
\end{gather*}
We can take the non-constant eigenvalue $\lambda_1(x)$ as one separation coordinate on~$S^3$.  However, since the other two
eigenvalues of~$L$ are equal, we can not recover all separation coordinates from the special conformal Killing tensor~$L$ in
this case.  We will have to consider another Killing tensor in the St\"ackel system instead~-- one with simple eigenvalues.

\looseness=-1
Before we do so, let us collect some facts from previous sections for the case $\Lambda_0\!<\!\Lambda_1\!=\!\Lambda_2\!=\!\Lambda_3$.
\begin{Proposition}
	Consider the following objects associated to a system of separation coordinates on $S^3$:
	\begin{itemize}\itemsep=0pt
		\item the St\"ackel system,
		\item a special conformal Killing tensor $L$ whose eigenvalues are constant on the coordinate hypersurfaces,
		\item the Benenti system $K(\lambda)=\Adj(L-\lambda g)$ of $L$,
		\item the special Killing tensor $K=L-(\tr L)g$ in the St\"ackel system,
		\item the St\"ackel line in the KS-variety,
		\item the trace free symmetric tensor $\hat L$ on $V$ which restricts to $L$,
		\item the vector $\vec n=(n_1,n_2,n_3)\in{\mathbb R}^3$ parametrising $\hat L$ via \eqref{eq:Lambda(n)}.
	\end{itemize}
	Then the following statements are equivalent:
	\begin{enumerate}\itemsep=0pt
		\item[$1.$] 
The St\"ackel system consists of extensions of Killing tensors from $S^2\subset S^3$
		                           $($up to the addition of a multiple of the metric$)$.
		\item[$2.$] 
$L$ has a multiple eigenvalue $($which is constant$)$.
		\item[$3.$] 
$K(\lambda)$ does not span the St\"ackel system.
		\item[$4.$] 
$K$ corresponds to a singular skew symmetric point in the KS-variety.
		\item[$5.$] 
The St\"ackel line lies on an isokernel plane in the KS-variety.
		\item[$6.$] 
$\hat L$ has a triple eigenvalue.
		\item[$7.$] 
$\lvert n_1\rvert=\lvert n_2\rvert=\lvert n_3\rvert$.
	\end{enumerate}
\end{Proposition}
\begin{proof}
	The equivalence of~6
and~7
is a direct consequence of~\eqref{eq:Lambda(n)}.  The equivalence to the
statements~1--5 follows from
Section~\ref{sec:extension},
Proposition~\ref{prop:eigenvalues},
Proposition~\ref{prop:triple},
Proposition~\ref{prop:singular}
	and
Lemma~\ref{lem:isokernel}, respectively.
\end{proof}

To make the extension of Killing tensors from $S^2\subset S^3$ more explicit, consider a Killing tensor on~$S^2$ with algebraic
curvature tensor~$R$ and Ricci tensor~$T$ and choose a basis in which $T$ is diagonal with diagonal elements~$t_1$,~$t_2$ and~$t_3$.  Then~$R$ is also diagonal and, if~$T$ is trace free, given by
\begin{gather*}
	R_{1212} =t_1+t_2=-t_3,\qquad
	R_{2323} =t_2+t_3=-t_1,\qquad
	R_{3131} =t_3+t_1=-t_2 .
\end{gather*}
In Section~\ref{sec:extension} we saw that this Killing tensor extends to a Killing tensor on $S^3$ corrsponding to the
point in the KS-variety that is def\/ined by the KS-matrix
\[
	V_t\coloneq
	\begin{pmatrix}
		t_3-t_2&-t_3&t_2\\
		t_3&t_1-t_3&-t_1\\
		-t_2&t_1&t_2-t_1\\
	\end{pmatrix}
	=
	t_1V_{+1}+t_2V_{+2}+t_3V_{+3} .
\]
If $T$ is not trace free, we obtain the same result.  This shows that the eigenvalues of the Ricci tensor of a Killing tensor on
$S^2$ induce barycentric coordinates on the isokernel planes in the KS-variety.

In Section~\ref{sec:S2} we have already computed the eigenvalues $\lambda_2(x)$ and $\lambda_3(x)$ of a Killing tensor on~$S^2$, parametrised by the eigenvalues of its Ricci tensor~$T$.  We see that if one perturbs the special Killing tensor on~$S^3$
within a St\"ackel system, the constant double eigenvalue $\lambda_2=\lambda_3$ splits up into two simple eigenvalues
$\lambda_2(x)$ and $\lambda_3(x)$, which come from separation coordinates on~$S^2$ under an embedding $S^2\subset S^3$.  On the
other hand, the eigenvalue $\lambda_1(x)$ remains unchanged.  Together, they def\/ine separation coordinates on~$S^3$.  We denote
this by writing the multiplicities of $\hat L$ as $(0(123))$ respectively $(0(1(23)))$, depending on whether the multiplicities
of the Ricci tensor $T$ are $(123)$ or $(1(23))$.  This yields the remaining two systems of separation coordinates:

{\bf (0(1(23))) Spherical coordinates.}
These are extensions of spherical coordinates on $S^2$ and correspond to St\"ackel lines through a face center and a vertex of the
octahedron in the KS-variety.

{\bf (0(123)) Lam\'e subgroup reduction.}
These are extensions of elliptic coordinates on $S^2$ and correspond to St\"ackel lines through a face center, but not a vertex of
the octahedron in the KS-variety.

\subsection{The classif\/ication space}

The non-generic separation coordinates are characterised by $n_1=\pm n_2$, $n_2=\pm n_3$ and $n_3=\pm n_1$.  This def\/ines six
projective lines that divide the projective plane of skew symmetric matrices in the KS-variety into twelve congruent triangles.
The four singular skew symmetric points are given by $\lvert n_1\rvert=\lvert n_2\rvert=\lvert n_3\rvert$ and each of the
triangles has two of them as vertices.  Upon a~blow-up in these points, the triangles become pentagons.  Any of these pentagons
constitutes a~fundamental domain for the $S_4$-action.  Regarding our notation for the multiplicities of~$\hat L$, we realise
that the vertices are labeled by the dif\/ferent possibilities to parenthesise the ordered set of eigenvalues of~$\hat L$ with
pairs of correctly matching parentheses and that the edges correspond to a~single application of the associativity rule.  This
naturally identif\/ies the classif\/ication space for separation coordinates on~$S^3$ with the associahedron~$K_4$, also known as
Stashef\/f polytope~\cite{Stasheff63}.

\begin{figure}[h]
	\centering
	\includegraphics[width=.55\textwidth]{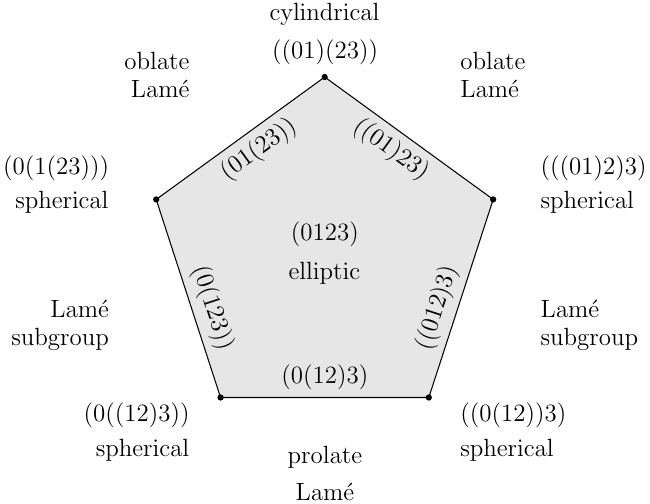}
	\caption{Identif\/ication of the moduli space of separation coordinates on $S^3$ with the associahed\-ron~$K_4$.}
	\label{fig:associahedron}
\end{figure}

\subsection*{Acknowledgements}

I would like to express my gratitude to Robert Milson for his motivation and the inspiring discussions about my f\/indings.  I~would also like to thank Alexander P.~Veselov for pointing out the link between my solution and moduli spaces of stable curves.
Finally, I would like to thank the anonymous referees, who helped to improve the paper considerably with their comments and
additional references.

\addcontentsline{toc}{section}{References}
\LastPageEnding

\end{document}